\documentclass[letterpaper, preprint, paper,11pt]{AAS}	

\usepackage{bm}
\usepackage{amsmath}
\usepackage{amssymb}
\usepackage{comment}
\usepackage{subcaption}
\usepackage[colorlinks=true, pdfstartview=FitV, linkcolor=black, citecolor= black, urlcolor= black]{hyperref}
\usepackage{overcite}
\usepackage{footnpag}			      	

\PaperNumber{26-931}

\begin{document}

\title{Transformer-Informed Trajectory Optimization for Relative Motion in Cislunar Orbits}

\author{Walter J. Manuel\thanks{PhD Candidate, Department of Aeronautics and Astronautics, Stanford University, Stanford, CA 94305.},
\ Yuji Takubo\footnotemark[1],
\ and Simone D'Amico\thanks{Associate Professor, Department of Aeronautics and Astronautics, Stanford University, Stanford, CA 94305.}
}

\maketitle{}

\begin{abstract}
Autonomous spacecraft guidance and control requires a fast solution to non-convex trajectory optimization, which can be accelerated by providing a near-optimal initial guess to an optimization protocol, i.e., warm-starting.
A robust warm starting method is especially useful for rendezvous, proximity operations, and docking (RPOD) in cislunar space, where the underlying dynamics become severely nonlinear and chaotic compared to those in Earth orbit, especially at perilune.
This paper extends the Autonomous Rendezvous Transformer (ART), a transformer-based warm-start trajectory generation method, to cislunar RPOD scenarios for the first time. 
To accurately and reliably solve the nonconvex optimal control problems (OCPs) posed by these scenarios, a new and enhanced version of ART, ART-TWIN (Two-Way INference), is introduced. 
Inspired by forward-backward shooting methods used in other trajectory design applications, ART-TWIN autoregressively generates two arcs, one from the initial state and one from the desired terminal state, that are patched together at the midpoint of the timeseries.
When evaluated on a set of simulated rendezvous scenarios that are initialized at perilune, ART-TWIN is demonstrated to substantially accelerate convergence and increase feasibility guarantees when used as a warm-start to sequential convex programming (SCP), compared to convex relaxations and the original ART.
These results illustrate the necessity of ART-TWIN's dual-arc generation to enable the viability of and gain benefits from using transformer-based warm-start methods in the most challenging areas of the cislunar dynamical regime.
\end{abstract}

\section{Introduction}

Mankind's long overdue return to the Moon is drawing closer, accompanied by both the promise of a flourishing cislunar space economy, and the deep technical challenges inherent in developing the orbital infrastructure to support such an ecosystem. 
As numerous nations and companies race to establish a sustained presence on the lunar surface, spearheaded by NASA's Artemis program and Moon Base initiative, assets in lunar orbit and the broader regime of cislunar space will serve essential roles in supporting crewed and robotic surface operations.
Among other needs, these roles include providing rendezvous, proximity operations, and docking (RPOD) services of both manned and unmanned spacecraft in cislunar space, and enabling the requisite space domain awareness (SDA) necessary to safely and securely execute lunar exploratory, scientific, and commercial objectives.
SDA and RPOD missions featuring distributed space systems (DSS), meaning systems involving two or more spacecraft, have become increasingly commonplace in Earth orbit, and thus it logically follows that advancements in DSS technology will soon be extended to cislunar space as well.
However, DSS operating in the cislunar regime are subject to the influence of both Earth and Moon gravity, where the governing dynamics become more nonlinear and chaotic.

This increase in dynamical complexity renders relative motion models formulated based on Keplerian dynamics, such as the Hill-Clohessy-Wiltshire equations\cite{clohessy_terminal_1960}, the Yamanaka-Ankersen state transition matrix (STM) \cite{yamanaka_new_2002}, and the Koenig-Guffanti-D'Amico STM \cite{koenig2017new}, less applicable \cite{hunter_optimal_2025}.
For this purpose, relative dynamics of the Circular Restricted Three-Body Problem (CR3BP), Elliptic Restricted Three-Body Problem (ER3BP), and Bi-Circular Restricted Four-Body Problem (BCR4BP) have been derived in the local-vertical local-horizontal (LVLH) frame \cite{franzini_relative_2019, khoury_relative_2022, khoury_relative_2024} and the velocity-normal-binormal (VNB) frame \cite{takubo_safe_2026, vela2026application}.

Various optimal guidance, navigation, and control (GNC) strategies are available to be adapted for use in cislunar space.
One such example leverages reachable set theory to solve convex optimal control problems (OCPs) via a computationally-efficient impulsive control algorithm. \cite{hunter_optimal_2025}.
However, different methods are required to address more challenging and realistic cislunar RPOD scenarios that involve non-convex objectives or constraints, such as keep-out-zone, approach ellipsoid, illumination, and field-of-view constraints \cite{guffanti_transformers_2024, {malyuta_thesis_2021}, khoury_relative_2024}.
The constraints involving passive safety are of particular importance to protect the RPOD-performing spacecraft \cite{breger2008safe}, as well as for the broader sustainability of the lunar environment.
The Earth's atmosphere provides a convenient avenue for orbital debris disposal, as objects in low Earth orbit (LEO) can naturally decay due to atmospheric drag and eventually burn up upon reentry. In contrast, since no atmosphere exists around the Moon, any debris generated from a cislunar collision could either stay indefinitely in orbit, posing a threat to other spacecraft, or could potentially strike the lunar surface, creating a cloud of sharp regolith particles and disrupting surface operations \cite{black_fragmentation_2025}.

Other constraints are sometimes self-imposed by mission designers to make the problem more tractable.
The difficulty of maneuvering at perilune (sometimes referred to as periselene or pericynthion), where spacecraft motion is fastest and the underlying multi-body dynamics are most chaotic, is often discussed in cislunar GNC literature, with many authors opting to simply limit maneuvers to apolune and near-apolune locations \cite{guzzetti_stationkeeping_2017,shimane_autonomous_2026}. 
This is a valid approach for some scenarios, such as long-term station keeping or reconfigurations \cite{foss_long-duration_2025}. 
However, such restrictions are incompatible with maintaining robust cislunar SDA.
The capability to rapidly react and respond on demand, including at perilune, for tasks such as in-space inspection, surface observation, or other needs, is necessary to ensure the safety of assets on the lunar surface and in orbit.

Combining these physics-driven and mission-driven factors together yields a dynamically challenging and non-convex OCP to solve.
While such non-convex OCPs are not new in spacecraft rendezvous and trajectory design, they are traditionally solved offline using computationally expensive methods \cite{malyuta2020fast, zhang2023stochastic, starek_real-time_2016}.
The proliferation of DSS requires improved onboard maneuver planning in order to provide more globally optimal solutions to the highly constrained, non-convex RPOD problems they are used for.
This higher level of autonomy becomes even more vital in cislunar space, given the communication latency between the Earth and Moon, as well as the communications blackout zone on the far side of the Moon\cite{shimane_autonomous_2026,liyanaarachchi_6g_2024}.
Recently, machine learning (ML) and artificial intelligence (AI) techniques have been applied to these autonomous spacecraft trajectory optimization problems to improve performance metrics such as fuel consumption and computational speed.\cite{guffanti_transformers_2024, banerjee_learning-based_2020, briden_constraint-informed_2025}.
Employing ML frameworks that are designed to learn and predict sequences, such as reinforcement learning or transformers, could improve the onboard computation of optimal trajectories for DSS operating in the cislunar regime.

The Autonomous Rendezvous Transformer (ART) takes an ML approach to solving the problem of trajectory optimization \cite{guffanti_transformers_2024}. 
Transformers are a class of neural network architectures which have gained much recognition for their use in natural language processing and other related applications \cite{vaswani_attention_2017}. 
Yet fundamentally, transformers are sequence predictors, making them applicable for trajectory generation. Rather than completely replacing traditional trajectory optimization methods, ART is trained to generate near-optimal trajectories, including both states and control inputs, that can be used as warm-starts for sequential convex programming (SCP) solvers. Having a more accurate warm-start can improve computational runtime, constraint satisfaction, and global optimality convergence when solving optimal control problems for RPOD scenarios.

ART was previously demonstrated in simulated LEO RPOD scenarios, where it compared favorably to other warm-start methods that used convex relaxations \cite{guffanti_transformers_2024, celestini_transformer-based_2024, takubo_towards_2025}. ART has also shown promise as a tool for rapid mission design\cite{takubo_agile_2026}.
However, up to this point, all usage of ART for spacecraft RPOD applications has been within the confines of the perturbed two-body problem, and has exclusively relied on linear dynamics models to characterize relative spacecraft motion.
Therefore, its ability to deal with the chaotic and highly nonlinear dynamics present in cislunar space has yet to be proven.
In order to fulfill a vision where ART is able to be fully deployed onboard spacecraft and used for autonomous operations in any domain, more efforts are needed to accelerate its convergence and robustness.
One option to achieve this goal is to incorporate a forward-backward shooting scheme, which has been widely successful in the field of low-thrust interplanetary trajectory optimization at reducing sensitivity, whether for direct methods \cite{vasile2003optimizing, takubo2023optimization}, indirect methods \cite{sidhoum_low_2023,sidhoum2024indirect}, or hybrid methods \cite{pierson1994three}.

With forward-backward shooting as a key enabler, this work seeks to make three main contributions that expand the scope and improve the performance of ART:
\begin{enumerate}
    \item A new version of ART, ART-TWIN (Two-Way INference), is developed by enhancing the original ART to include forward-backward dual-arc trajectory generation, which makes it possible for ART to achieve reliable convergence in dynamically chaotic regimes.
    \item ART-TWIN is applied to cislunar RPOD scenarios using CR3BP nonlinear relative dynamics models. 
    This marks the first extension of ART to applications beyond Earth orbit, as well as the first time ART has used a nonlinear dynamics model for spacecraft trajectory generation. 
    Furthermore, ART-TWIN is shown to be capable of usage at perilune, where the dynamics are most nonlinear and challenging. 
    \item In batch testing of simulated cislunar rendezvous scenarios, ART-TWIN demonstrates comparable optimality performance and significant efficiency gains, notably a 1.59x faster median runtime, when compared to a convex relaxation warm-start. When applied to the same test cases, the warm-start provided by the original version of ART (with forward-shooting only) fails to converge via SCP in over 60\% of instances.
\end{enumerate}
The remainder of this paper is organized as follows: First, relevant background details on the CR3BP and ART are given. Then, the specific methodology featured by this work is introduced, including the forward-backward shooting architecture of ART-TWIN, the cislunar RPOD OCP formulation, and the incorporation of SCP. Finally, results and analysis obtained from applying ART-TWIN to a simulated cislunar rendezvous scenario are presented, followed by conclusions and proposed directions for future work.

\section{Background}

\subsection{Cislunar Relative Dynamics}

\subsubsection{Circular Restricted Three-Body Problem}

This paper considers a system comprised of two spacecraft, a chief and a deputy, operating in the cislunar system and influenced by the gravitational forces of both the Earth and the Moon. 
As a result, the CR3BP is employed as a reduced-order cislunar dynamics model to characterize the absolute and relative motion of the two spacecraft.
In the CR3BP, a spacecraft of relatively negligible mass is accompanied by two bodies of significant mass, known as primaries, such as the Earth and the Moon. 
Under standard CR3BP simplifying assumptions, the two larger bodies are modeled as point masses and move in circular orbits about their common center of mass. 
It is convenient and conventional in CR3BP dynamic analysis to nondimensionalize the units of distance and time, and to express motion in a synodic rotating frame that has the same inertial angular velocity as the two primaries\cite{howell1984three}.

\begin{figure}[htb!]
    \centering
    \begin{subfigure}[b]{.3\textwidth}
        \centering
        \includegraphics[width=\textwidth]{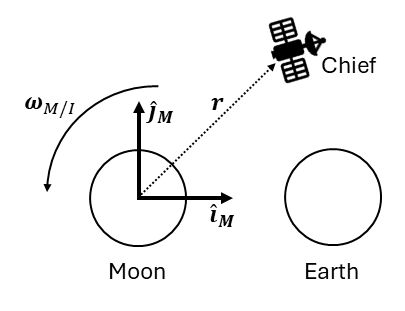}
        \captionsetup{width=1.5\linewidth, justification=centering}
        \caption{Moon-centered synodic reference frame}
        \label{fig:moon_frame}
    \end{subfigure}\hspace{0.3\textwidth}
    \begin{subfigure}[b]{.2\textwidth}
        \centering
        \includegraphics[width=0.9\textwidth]{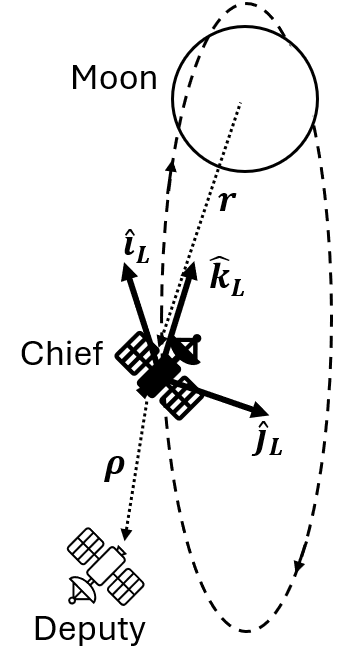}
        \captionsetup{width=1.5\linewidth, justification=centering}
        \caption{LVLH reference frame}
        \label{fig:lvlh_frame}
    \end{subfigure}  
    \caption{Key reference frames, vectors, and quantities used in this work. The chief spacecraft is depicted as being on a Near-Rectilinear Halo Orbit (NRHO).\cite{hunter_optimal_2025}}
    \label{fig:frames} 
\end{figure}

The system of deputy and chief spacecraft is visualized in Figure \ref{fig:frames}, along with the two key reference frames employed in this paper. 
The chief spacecraft is depicted as being on a Near-Rectilinear Halo Orbit (NRHO). 
NRHOs are a subclass of halo orbits, three-dimensional periodic orbits that naturally form about the equilibrium points of the CR3BP system \cite{zimovan2017characteristics}. 
%
The Moon-centered synodic rotating frame $\mathcal{M}$, as shown in Figure \ref{fig:moon_frame}, is defined by unit vectors $\hat{\bm{i}}_{M}$, $\hat{\bm{j}}_{M}$, and $\hat{\bm{k}}_{M}$, where $\hat{\bm{i}}_{M}$ points from the Moon to the Earth, $\hat{\bm{k}}_{M}$ is aligned with the angular momentum vector of the Moon with respect to the Earth, and $\hat{\bm{j}}_{M}$ completes the right-handed triad. 
An LVLH frame $\mathcal{L}$, centered on the chief spacecraft and anchored to the Moon, is also employed as shown in Figure \ref{fig:lvlh_frame}. Frame $\mathcal{L}$ is defined by unit vectors $\hat{\bm{i}}_{L}$, $\hat{\bm{j}}_{L}$, and $\hat{\bm{k}}_{L}$, where $\hat{\bm{k}}_{L}$ points from the chief to the Moon, $\hat{\bm{j}}_{L}$ is anti-parallel to the angular momentum vector of the chief with respect to the Moon, and $\hat{\bm{i}}_{L}$ completes the right-handed triad.
The LVLH frame is sometimes alternatively defined as the Radial/Tangential/Normal (RTN) frame, using a permutation of the same unit vectors, $\{ -\hat{\bm{k}}_L, \hat{\bm{i}}_L, -\hat{\bm{j}}_L \}$.
The absolute position of the chief spacecraft relative to the Moon, depicted in both Figures \ref{fig:moon_frame} and \ref{fig:lvlh_frame}, is denoted as $\bm{r}$, and the relative position of the deputy spacecraft with respect to the chief, depicted in Figure \ref{fig:lvlh_frame}, is expressed as $\bm{\rho}$.

\subsubsection{Nonlinear Cartesian Relative Dynamics}
The cislunar relative dynamics model derived by Franzini et al. uses an LVLH frame anchored to the Moon\cite{franzini_relative_2019}.
The nonlinear form of these equations of motion is given as 
\begin{align} \label{eq:nonlinear_cr3bp_rel_dyn}
    [\bm{\ddot{\rho}}]^{\mathcal{L}} = 
    &-2\bm{\Omega}_{L/I}[\bm{\dot{\rho}}]^{\mathcal{L}}
    - \left([\bm{\dot{\Omega}}_{L/I}]^{\mathcal{L}} + \bm{\Omega}^{2}_{L/I}\right)\bm{\rho}
    + \mu \left(\frac{\bm{r}}{\| \bm{r} \|^3} - \frac{\bm{r + \rho}}{{\| \bm{r + \rho} \|}^3}\right)\\
    &+ \left( 1 - \mu \right) \left(\frac{\bm{r} + \bm{r}_{em}}{\| \bm{r} + \bm{r}_{em} \|^3} - \frac{\bm{r + \rho} + \bm{r}_{em}}{{\| \bm{r + \rho} + \bm{r}_{em} \|}^3}\right), \notag
\end{align}
where $\bm{r}_{em}$ is the position vector of the Earth relative to the Moon, $\bm{r}$ is the position vector of the chief spacecraft relative to the Moon, and $\mu$ is the mass parameter of the CR3BP system \cite{howell1984three}.
The skew-symmetric matrices $\bm{\Omega}_{L/I} \in \mathbb{R}^{3\times3}$ and $[\bm{\dot{\Omega}}_{L/I}]^{\mathcal{L}}\in \mathbb{R}^{3\times3}$ are associated with the angular velocity vector of frame $\mathcal{L}$ relative to the inertial frame $\mathcal{I}$, $\bm{\omega}_{L/I}$, and the angular acceleration vector $[\bm{\dot{\omega}}_{L/I}]^{\mathcal{L}}$, respectively. 
These two matrices are defined as
\begin{align}\label{skew_matrix}
    \bm{\Omega}_{L/I} = \begin{bmatrix}
        \quad0 & -\omega^{z}_{L/I} & \quad\omega^{y}_{L/I}\\
        \quad\omega^{z}_{L/I} & \quad0 & -\omega^{x}_{L/I}\\
        -\omega^{y}_{L/I} & \quad\omega^{x}_{L/I} & \quad0
    \end{bmatrix},\quad
    [\bm{\dot{\Omega}}_{L/I}]^{\mathcal{L}} = \begin{bmatrix}
        \quad0 & -\dot{\omega}^{z}_{L/I} & \quad\dot{\omega}^{y}_{L/I}\\
        \quad\dot{\omega}^{z}_{L/I} & \quad0 & -\dot{\omega}^{x}_{L/I}\\
        -\dot{\omega}^{y}_{L/I} & \quad\dot{\omega}^{x}_{L/I} & \quad0
    \end{bmatrix}.
\end{align}
\subsubsection{Linearized Cartesian Relative Dynamics}
By linearizing the gravitational acceleration terms in Eq. (\ref{eq:nonlinear_cr3bp_rel_dyn}) about the chief's position in the LVLH frame, the free relative dynamics in the CR3BP system can be formulated as a linear time-variant (LTV) system of equations in state space form as $\bm{\dot{x}} = \bm{A}(t)\bm{x}$, where the state vector $\bm{x}$ is the Cartesian position and velocity of the deputy spacecraft in the LVLH frame, $\begin{bmatrix} \bm{\rho}\\ \bm{\dot{\rho}} \end{bmatrix}$. 
The plant matrix $\bm{A}(t)$ is defined as
\begin{align}
    \bm{A}(t) = 
    \begin{bmatrix}
        \bm{0}_{3}                &        \bm{I}_3 \\
        \bm{A_{\dot{\rho}\rho}}(t)  &  -2\bm{\Omega}_{L/I}(t)
    \end{bmatrix}. \label{eq: cr3bp A}
\end{align}
where $\bm{0}_3$ and $\bm{I}_3$ are $3\times3$ zero and identity matrices.
The $3\times3$ matrix $\bm{A_{\dot{\rho}\rho}}$ concisely combines the terms in the lower left quadrant, defined as
\begin{align} \label{eq:linear_cr3bp_rel_dyn}
    \bm{A_{\dot{\rho}\rho}}(t) = &-[\bm{\dot{\Omega}}_{L/I}]^{\mathcal{L}} -\bm{\Omega}^{2}_{L/I} 
    - \frac{\mu}{\| \bm{r} \|^3}\left(\bm{I} - 3\frac{\bm{r}\bm{r}^\top}{{\| \bm{r} \|}^2}\right)\\
    &- \frac{1 - \mu}{\| \bm{r} + \bm{r}_{em} \|^3}\left(\bm{I} - 3\frac{(\bm{r} + \bm{r}_{em})(\bm{r} + \bm{r}_{em})^\top}{{\| \bm{r} + \bm{r}_{em} \|}^2}\right), \notag
\end{align}
%
The linearized model given in Eqs. (\ref{eq: cr3bp A}) and (\ref{eq:linear_cr3bp_rel_dyn}) is expected to perform within an acceptable range of error at apolune, but not at perilune, based on the validation performed by its original authors \cite{franzini_relative_2019}.

\subsection{Autonomous Rendezvous Transformer}
ART is a transformer-based trajectory generation algorithm. 
It was developed to leverage the generation capability of the causal transformer to autoregressively predict and then generate each step of a trajectory. 
A key enabler of ART is the tokenization of trajectory information, which allows sequences of states, control actions, and performance metrics to be encoded together and modeled in a way that is compatible with the transformer.
More specifically, a tokenized ART trajectory with $N$ discretized time steps over the horizon is represented as 
\begin{align} \label{eq:tau}
    \bm{\tau}_{1:N} = \left\{ \bm{x}_1, \bm{u}_1, r_1, c_1, ..., \bm{x}_N, \bm{u}_N, r_N, c_N \right\},
\end{align}
where at time step $k$, $\bm{x}_k \in \mathbb{R}^{s}$ is the state, $\bm{u}_k \in \mathbb{R}^{a}$ is the control action, and $r_k$ and $c_k$ are two key performance metrics that give users control over trajectory design.
The reward-to-go (RTG) $r_k \in \mathbb{R}$ is defined as the cost incurred from the current epoch $t_k$ to the terminal epoch $t_N$:
\begin{align}
    r_k = -\sum^{N}_{j=k} J(\bm{x}_j,\bm{u}_j).
\end{align}
The generic cost function $J$ could represent any metric that is desired to be minimized along the trajectory, such as fuel consumption. Given that $\bm{\mathcal{C}}_k$ is the feasible set of the states and actions, the constraint-to-go (CTG) $c_k \in \mathbb{N}$ is defined as the sum of steps along the discretized trajectory that contain constraint violations, from $t_k$ to $t_N$:
\begin{align}\label{eq:rtg_c}
    c_k & = \sum_{j=k}^{N} {\textsf{C}}_j, \\
    {\textsf{C}}_j &= 
    \begin{cases}
        1 & \text{if  } {\exists \bm{x}_j, \bm{u}_j} \notin \bm{\mathcal{C}}_k \\
        0 & \text{otherwise}
    \end{cases}.
\end{align}
The type of problem that ART is built to help solve is a fixed-time, non-convex OCP, which is formulated generically as follows:
\begin{align} \label{eq:OCP}
& \underset{\bm{x}_k, \bm{u}_k}{\text{minimize}} & & \sum_{k}^{N} {J}(\bm{x}_k, \bm{u}_k) &\\ \nonumber
& \text{subject to} & & \bm{x}_{k+1} = \bm{f}\left(\bm{x}_{k}, \bm{u}_{k} \right)  & \forall k \in [1, N-1] \\ \nonumber
&&& (\bm{x}_k, \bm{u}_k) \in \bm{\mathcal{C}}_k & \forall k \in [1, N] \\ \nonumber
&&& \bm{x}_1 = \bm{x}_i, \ \bm{x}_N = \bm{x}_f, 
\end{align}
where $\bm{f}(\bm{x}_k, \bm{u}_k): \mathbb{R}^{s+a} \rightarrow \mathbb{R}^{s}$ represents the system dynamics (which can be linear or nonlinear).
In this problem setup, both initial and terminal conditions are fixed to pre-defined values, given as $\bm{x}_i$ and $\bm{x}_f$, respectively.

\subsubsection{Dataset Generation}

The first step in the implementation pipeline of ART is to generate a dataset suitable for effective transformer training and testing.
The dataset consists of $N_d$ trajectories obtained by repeatedly solving diverse instances of the optimization problem in Eq. (\ref{eq:OCP}) via SCP, with each solved trajectory sample containing both state and action timeseries.
To lend practicality to the setup, the parameters of the OCP can be set based on a realistic designed mission scenario, with some boundary conditions, such as the initial state vector of the deputy spacecraft, randomized to provide the aforementioned diversity to the dataset.
By exposing the transformer model to optimal or near-optimal trajectories with varying initial conditions and performance metrics, it is able to learn the underlying dynamics and variations in trajectory behavior. This enables it to generate control policies at inference for boundary conditions outside the training set, while still maintaining robust performance.

\subsubsection{Model Training}

ART is trained via the standard teacher-forcing procedure commonly used in training sequence models.\cite{guffanti_transformers_2024, celestini_transformer-based_2024}
In this process, the loss function
\begin{equation}
\begin{aligned}
    {\mathcal{L}}(\tau) = \sum_{m=1}^{M} \sum_{k=1}^{N} \left( \Vert \bm{x}_{m,k} - \hat{\bm{x}}_{m,k} \Vert_2^2 + \Vert \bm{u}_{m,k} - \hat{\bm{u}}_{m,k} \Vert_2^2 \right)
\end{aligned},
\end{equation}
is minimized over a batch of $M$ trajectories, where $\hat{\bm{x}}$ and $\hat{\bm{u}}$ represent the predicted state and control action inferred by the transformer, respectively.

\subsubsection{Model Inference}

In the original architecture\cite{guffanti_transformers_2024}, ART is given three inputs when generating a trajectory: the initial state $\bm{x}_i$, initial reward-to-go $r_1$, and initial constraint-to-go $c_1$.
ART then infers the next control input $\bm{u}_k$, which is fed to a dynamics model  $\bm{f}(\bm{x}_k, \bm{u}_k)$ to give the next state $\bm{x}_{k+1}$. Then the new reward-to-go $r_k$, and new constraint-to-go $c_k$, are computed and fed back into ART along with the new state, and the inference process repeats. This autoregressive process continues until all $N$ states and control inputs in the discretized trajectory are computed. 
This full trajectory generated by ART is then used as an initial guess to warm-start SCP, which produces an optimal final solution to the trajectory optimization problem.
The final solution produced via SCP enforces feasibility and hard constraint satisfaction.
This description applies to the dynamics-informed mode of ART, which is the only mode applied in this work. A transformer-only mode in which ART infers successive states directly without any dynamical knowledge, albeit with significantly degraded performance compared to the dynamics-informed mode, was also developed in the original implementation of ART \cite{guffanti_transformers_2024}.

\section{Methodology}

\subsection{Forward-Backward Dual-Arc Generation}

\begin{figure}[htb!]
    \centering
    \includegraphics[width=0.6\linewidth]{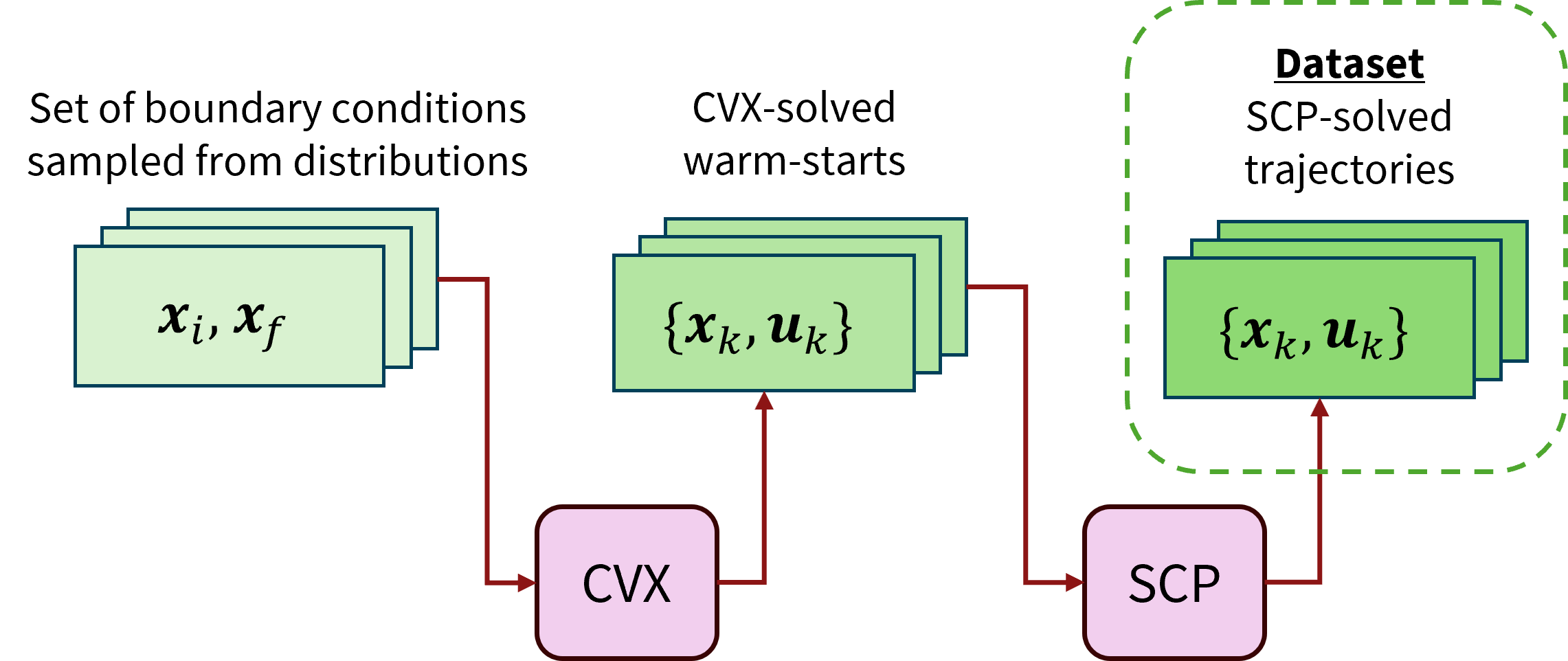}
    \caption{Flowchart of ART-TWIN dataset generation process. In the CVX block, the convex relaxation of the non-convex OCP is solved to provide a warm-start for the SCP-driven trajectory generation.}
    \label{fig:art_datagen}
\end{figure}

\begin{figure}[htb!]
    \centering
    \includegraphics[width=0.9\linewidth]{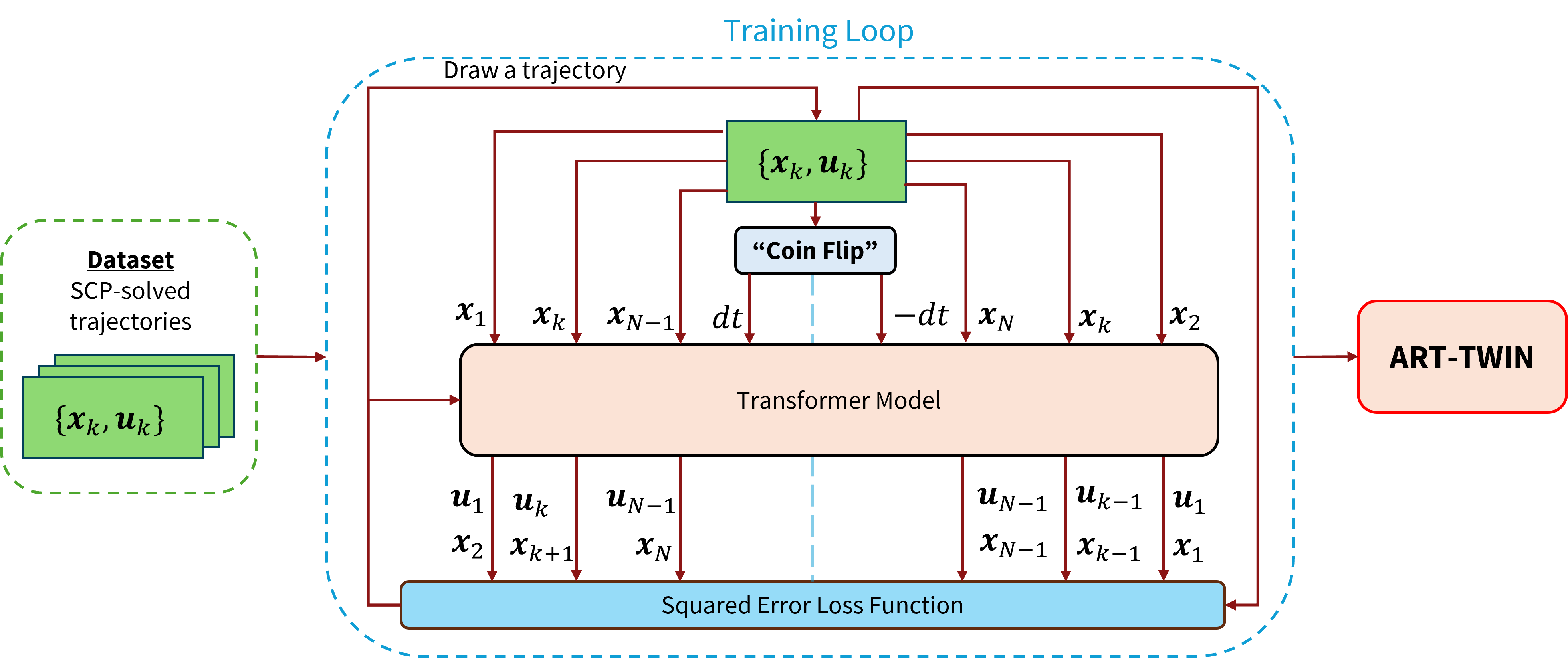}
    \caption{Flowchart of (dynamics-informed) ART-TWIN training process. At each step in the training loop, a new trajectory is randomly drawn with replacement from the training set. Based on the outcome of a 50-50 random draw (a "coin flip"), the selected trajectory is either assigned token $dt$ which routes the transformer model to train via forwards propagation starting from the initial state, or assigned token $-dt$ which routes the transformer model to train via backwards propagation starting from the terminal state.}
    \label{fig:art_training}
\end{figure}

\begin{figure}[htb!]
    \centering
    \includegraphics[width=1\linewidth]{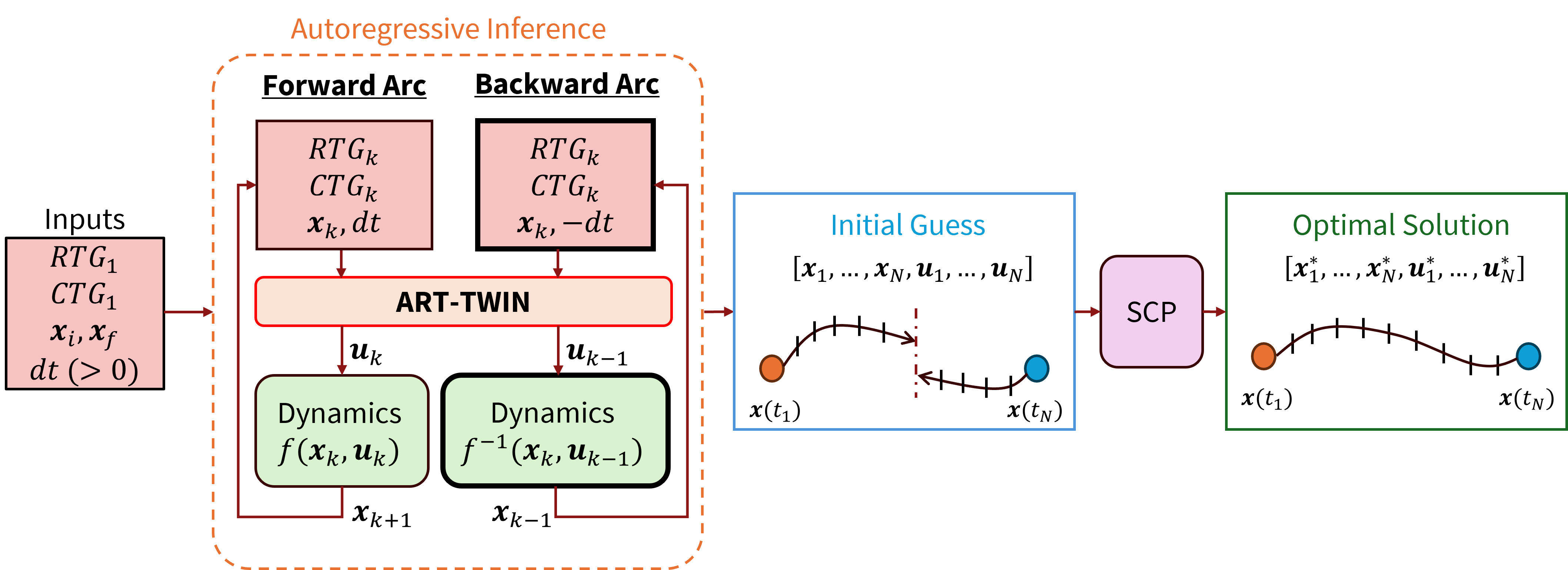}
    \caption{Flowchart of (dynamics-informed) ART-TWIN inference and warm-start pipeline.}
    \label{fig:art_workflow}
\end{figure}

The proposed method, ART-TWIN, integrates a forward-backward shooting scheme into ART.
The dataset generation of ART-TWIN is depicted in Figure \ref{fig:art_datagen}, its training process is shown in Figure \ref{fig:art_training}, and its inference workflow is illustrated in Figure \ref{fig:art_workflow}. 
Leveraging the fact that both boundary states are given, the autoregressive trajectory generation loop then propagates two arcs at inference: a forward arc, starting from $\bm{x}_1$ and propagated along the forward dynamics with increment $dt$, and a backward arc, starting from $\bm{x}_N$ and propagated along the backward dynamics with increment $-dt$. The backwards dynamics are defined as $\bm{x}_{k-1} = \bm{f}^{-1}(\bm{x}_k, \bm{u}_{k-1})$, where in practice, $\bm{f}^{-1}$ simply denotes integrating the discrete nonlinear dynamics in reverse, from timestep $t_k$ to $t_{k-1}$, as opposed to forward integration from $t_k$ to $t_{k+1}$.
Each arc is propagated until reaching a predefined common breakpoint, which is set in this work to be the midpoint of the rendezvous time horizon.
The combined trajectory of the two arcs is then patched and input to the SCP solver as an initial guess.

\begin{figure}[htb!]
    \centering
    \includegraphics[width=0.7\linewidth]{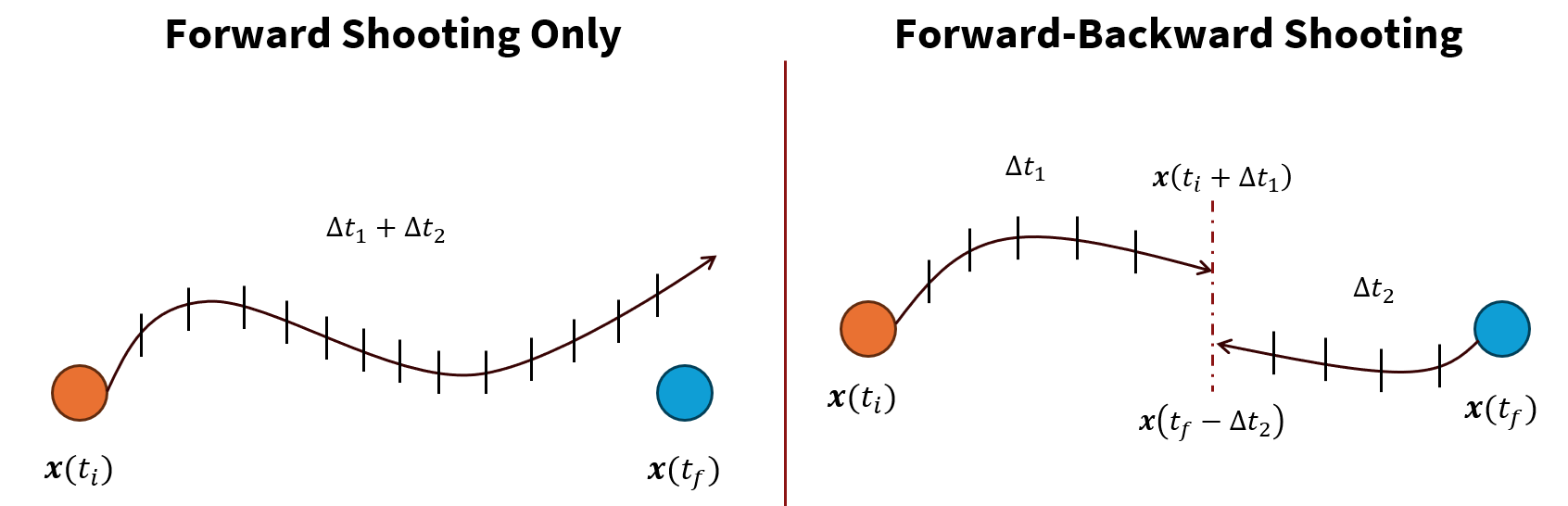}
    \caption{Forward vs. forward-backward shooting schemes for ART inference. ART-TWIN employs forward-backward inference.}
    \label{fig:forward_backward}
\end{figure}

The difference between the generated warm-start trajectories is illustrated in Figure \ref{fig:forward_backward}. 
Since the original ART generates the whole trajectory in forward time direction, the deviation accumulates over time and the terminal state often has a deviation from the desired state $\bm{x}_f$. 
This is especially salient for the chaotic dynamics of the CR3BP, where a small perturbation to the initial conditions can lead to a large variation in the terminal state.
In contrast, ART-TWIN provides two arcs with a breakpoint in the middle. 
Since the length of each generated arc is halved, and both boundary states are fixed, the warm-start is expected to have less deviation, and be less sensitive to initial conditions, thus improving the convergence of the following optimization process. 
This does not mean that the ART-TWIN warm-start has zero deviation, as there is still a defect between the two arcs, i.e., the states $x_{\frac{N}{2}}$ and $x_{\frac{N}{2} + 1}$ do not conform to the proper dynamical flow.
However, the presumed advantage of forward-backward shooting is that this defect is easier to patch and resolve through SCP iterations, as opposed to resolving terminal state error.
In order to realize the inference in both time directions, an additional argument $dt$ is added as an input token of ART-TWIN. 
While the original ART passes only reward-to-go ($r_k$), constraint-to-go ($c_k$), and the current state $\bm{x}_k$ at each time step, ART-TWIN receives explicit information of the propagation time. 
This extends the current capability of ART to varying time discretizations, and also enables inference to be made backward in time by training a single inference model.

To be clear on notation going forward, ART will describe the broad class of ART model versions, ART-0 will refer to the original version of ART (forward shooting only), and ART-TWIN will refer to the newly developed and implemented version of ART (with forward-backward dual-arc generation).

\subsection{Applying ART-TWIN to Cislunar Rendezvous}

The initial application and evaluation of the cislunar-capable ART-TWIN follows a similar pipeline to that of ART-0. 
First, ART-TWIN is trained on a batch of SCP solutions to non-convex CR3BP reconfiguration OCPs with a passive safety constraint \cite{guffanti2023passively}.
Then, the performance of ART-TWIN as a warm start is tested against a warm start trajectory obtained by solving the convex rendezvous problem which does not enforce the keep-out-zone constraint.
ART-TWIN is expected to equal or outperform the convex benchmark in terminal state accuracy, number and severity of constraint violations, control cost, and runtime.
ART-TWIN is also evaluated in comparison to ART-0 to demonstrate its improvements in convergence and efficiency.
In this training and testing pipeline, two OCP formulations are considered.

\subsubsection{Problem 1: Convex, Linear Rendezvous Problem}

The first OCP of interest is a convex relaxation of a fuel-optimal rendezvous problem, and is given as
\begin{align} \label{eq: OCP_linear}
& \underset{\bm{x}_k, \bm{u}_k}{\text{minimize}} & & \sum_{k=1}^{N} \| \bm{u}_k \|_2 &\\ \nonumber
& \text{subject to} & & \bm{x}_{k+1} = \bm{\Phi}(t_{k+1}, {t_k}) \bm{x}_k + \bm{B}_k\bm{u}_k  & \forall i \in [1, N-1] \\ \nonumber
&&& \bm{x}_1 = \bm{x}_i, \ \bm{x}_N = \bm{x}_f, \\ \nonumber
\end{align}
where $\bm{B}_k$ is the control input matrix $\begin{bmatrix} \bm{0}_3 \\ \bm{I}_3 \end{bmatrix}$, and $\bm{\Phi}(t_{k+1}, {t_k})$ is the STM from timesteps $t_{k}$ to $t_{k+1}$ created by taking the matrix exponential of Eq. (\ref{eq: cr3bp A}) computed by expanding the power series of $e^A$ to fifth order.
The solution to Eq. \eqref{eq: OCP_linear} is used as the warm-start to generate the training dataset of solved SCP trajectories, and as the baseline comparison to the warm-start provided by ART-TWIN.

\subsubsection{Problem 2: Non-Convex, Nonlinear Rendezvous Problem}

The second, more complex OCP considered in this work provides the basis for the trajectories used to train ART-TWIN, and requires the usage of an optimization protocol such as SCP, in order to iteratively solve it.
This nonlinear and nonconvex OCP is given as
\begin{align} \label{eq: OCP_nonlinear}
& \underset{\bm{x}_k, \bm{u}_k}{\text{minimize}} & & \sum_{k=1}^{N} \| \bm{u}_k \|_2 &\\ \nonumber
& \text{subject to} & & \bm{x}_{k+1} = \bm{f}\left(\bm{x}_{k}, \bm{u}_{k} \right)  + \bm{B}_k\bm{u}_k & \forall k \in [1, N-1] \\ \nonumber
&&& \bm{x}_{kj}^\top \bm{E}_{ae}^\top \bm{E}_{ae} \bm{x}_{kj} \geq 1 & \forall k \in [1, n_{gate}],\ \forall j \in [1, n_{safe}] \\ \nonumber
&&& \bm{x}_{kj}^\top \bm{E}_{koz}^\top \bm{E}_{koz} \bm{x}_{kj} \geq 1 & \forall k \in [n_{gate}+1, N],\ \forall j \in [1, n_{safe}] \\ \nonumber 
&&& \bm{x}_1 = \bm{x}_i, \ \bm{x}_N = \bm{x}_f \\ \nonumber
& \text{where} & & \bm{x}_{kj} = \left( \prod_{m=1}^{j-1} \bm{\Phi}(t_{m+1}, {t_m}) \right) \bm{x}_k & \forall j \in [1, n_{safe}].
\end{align}
The dynamics function $\bm{f}\left(\bm{x}_{k}, \bm{u}_{k} \right)$ is the discrete-time, nonlinear form of the cislunar relative equations of motion found in Eq. (\ref{eq:nonlinear_cr3bp_rel_dyn}), which is used as the simulation ground truth as well.
The matrices $\bm{E}_{ae} = \text{diag}([1/r_{ae}, 1/r_{ae}, 1/r_{ae},0,0,0])$ and $\bm{E}_{koz} = \text{diag}([1/r_{koz}, 1/r_{koz}, 1/r_{koz},0,0,0])$ denote the approach and keep-out-zone ellipsoids, with radii $r_{ae}$ and $r_{koz}$, respectively. 
These ellipsoids are used to formulate nonconvex passive safety constraints, which enforce free-drift passive safety, i.e., keeping the spacecraft from intersecting with the ellipsoid for $n_{safe}$ steps within a predefined safety horizon. 
The approach sphere constraint is active until step $n_{gate}$ of the rendezvous horizon time discretization, at which point the keep-out-zone becomes active for the remaining duration of the rendezvous horizon.

\subsection{Sequential Convex Programming}

SCP is a technique to solve non-convex OCPs by iteratively solving convexified subproblems until convergence to an optimal solution is achieved. 
This work uses the SCvx* algorithm as the backbone of its SCP solver \cite{Oguri_2023}.
The specific SCvx* hyperparameters used are provided in Table \ref{tab:scp_params} in the Appendix. 
SCvx* incorporates violations of the linearized nonconvex constraints into a penalty function in order to formulate an augmented Lagrangian. This guarantees the feasibility of convexified subproblems. Through each iteration of SCP, the trust region and penalty weights are both updated until the solver converges to a solution.
During the SCP routine, numerical stability is maintained by normalizing all variables with respect to the initial guess, with each element scaled to the range [-1,1].
Although theoretical convergence guarantees are discussed in various SCP algorithms, non-convergence is still a possibility, due to either lack of feasibility of the OCP, numerical instabilities, or exceeding the maximum allowed number of iterations.
When solving non-convex problems, it is possible for SCP to converge to different local optima based on the provided initial guess.
Hence, the input of a ``better" warm-start is critical in guiding SCP to converge to a more optimal solution. 

At each SCP iteration, the following convex subproblem is solved: 
\begin{align} \label{eq:cvx_subproblem}
    & \underset{\bm{x}_k, \bm{u}_k, \xi_k}{\text{minimize}} & & \sum_{k = 1}^{N} || \bm{u}_k ||_2  + \mathcal{J}_{\text{pen}}(\bm{\xi})& \\ \nonumber 
    & \text{subject to} & &  \bm{x}_{k+1} = f(\bar{\bm{x}}_k, \bm{u}_k) + \left.\frac{\partial \bm{f}}{\partial \bm{x}}\right|_{\bar{\bm{x}}_k} (\bm{x}_k - \bar{\bm{x}}_k) + \bm{B}_k\bm{u}_k  & \forall k \in [1, N) \\ \nonumber 
    &&& - \bm{a}^\top_k \bm{x}_k + b_k \leq \xi_k  & \forall k \in [1, N] \\ \nonumber 
    &&& \bm{x}_1 = \bm{x}_i, \ \bm{x}_N = \bm{x}_f, \\ \nonumber
    &&& \| \bm{x}_k - \bar{\bm{x}}_k\|_\infty \leq \epsilon &  \forall k \in [1, N] \\ \nonumber 
    &&& \| \bm{u}_k- \bar{\bm{u}}_k \|_\infty \leq \epsilon &  \forall k \in [1, N-1] \\ \nonumber 
    & \text{where}& & \bm{a}_{kj} = \bar{\bm{x}}_k^\top \bm{\Psi}^\top_{j} E_{ae} \bm{\Psi}_j & \forall k \in [1, n_{gate}+],\ \forall j \in [1, n_{safe}]  \\ \nonumber 
    &&& \bm{a}_{kj} = \bar{\bm{x}}_k^\top \bm{\Psi}^\top_{j} E_{koz} \bm{\Psi}_j & \forall k \in [n_{gate}+1, N],\ \forall j \in [1, n_{safe}]  \\ \nonumber 
    &&& b_{kj} = \sqrt{\bm{a}^\top_{kj} \bar{\bm{x}}_k}& \forall k \in [1, N],\ \forall j \in [1, n_{safe}]  \\ \nonumber
    &&& j^*_k = \underset{j}{\arg\min} \ b_{kj} & \forall k \in [1, N]\\ \nonumber
    &&& b_k \triangleq b_{kj^*} & \forall k \in [1, N]  \\ \nonumber
    &&& \bm{a}_k \triangleq \bm{a}_{kj^*} & \forall k \in [1, N]  \\ \nonumber
    &&& \bm{\Psi}_j = \prod_{m=1}^{j-1} \bm{\Phi}(t_{m+1}, {t_m}) & \forall j \in [1, n_{safe}].
\end{align}
The passive safety constraint is linearized by using the deputy's states from the previous SCP iteration, $\bar{\bm{x}}_k$, as a reference. It is important to note that nonlinear dynamics are still incorporated in this convexified subproblem. As shown in the second line of Eq. \eqref{eq:cvx_subproblem}, at each discretized step of the trajectory, the reference state is propagated forward to the next step using the nonlinear equations of motion, and a first-order Taylor series expansion with respect to the deputy's current state is taken, leaving the decision variables free for the solver to control. This ensures that the final converged solution will still be identical to a nonlinear ground truth propagation, within machine error tolerances. Further details on the slack variable $\bm{\xi}$, and the penalty function $\mathcal{J}_{\text{pen}}(\bm{\xi})$, can be found in the original documentation of SCvx* \cite{Oguri_2023}.

It is worth noting that it is possible for a dual- or multi-arc generation scheme to also be applied to the SCP process, either at the dataset generation step, inference step, or both.
While that is a viable avenue of investigation, this work focuses on gaining benefits by using ART to boost traditional optimization methods, rather than the converse.

\section{Results and Analysis}

This section details the experiments conducted to assess the performance of ART-TWIN in the context of its application to cislunar RPOD scenarios.\footnote{The experiments are conducted using an AMD Ryzen 9 7950X CPU, and an NVIDIA GeForce RTX 4090 GPU.}
\subsection{Cislunar Rendezvous Scenario}

\begin{figure}[htb!]
    \centering
    \begin{subfigure}[b]{0.48\textwidth}
        \centering
        \includegraphics[width=\linewidth]{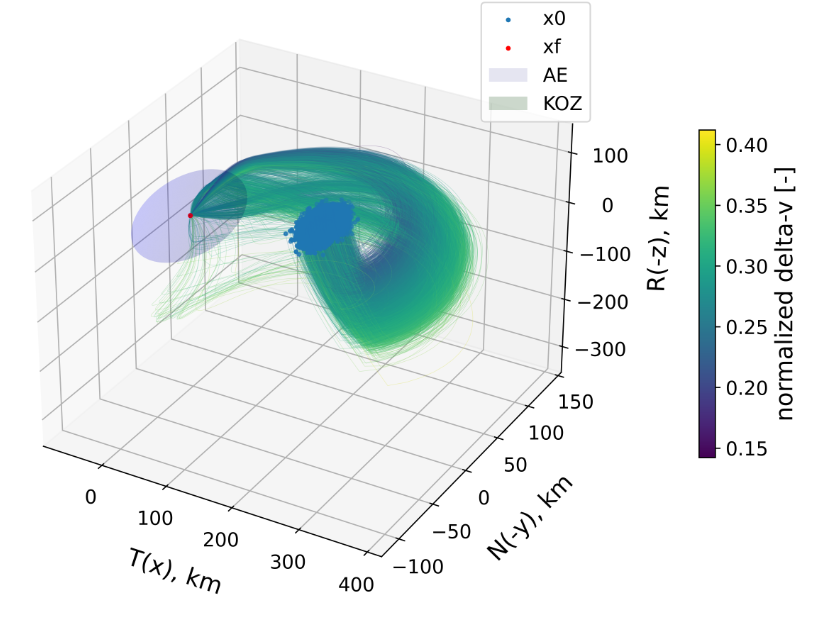}
        \caption{3D Plot, LVLH Frame}
        \label{fig:sub1}
    \end{subfigure}
    \hfill
    \begin{subfigure}[b]{0.48\textwidth}
        \centering
        \includegraphics[width=\linewidth]{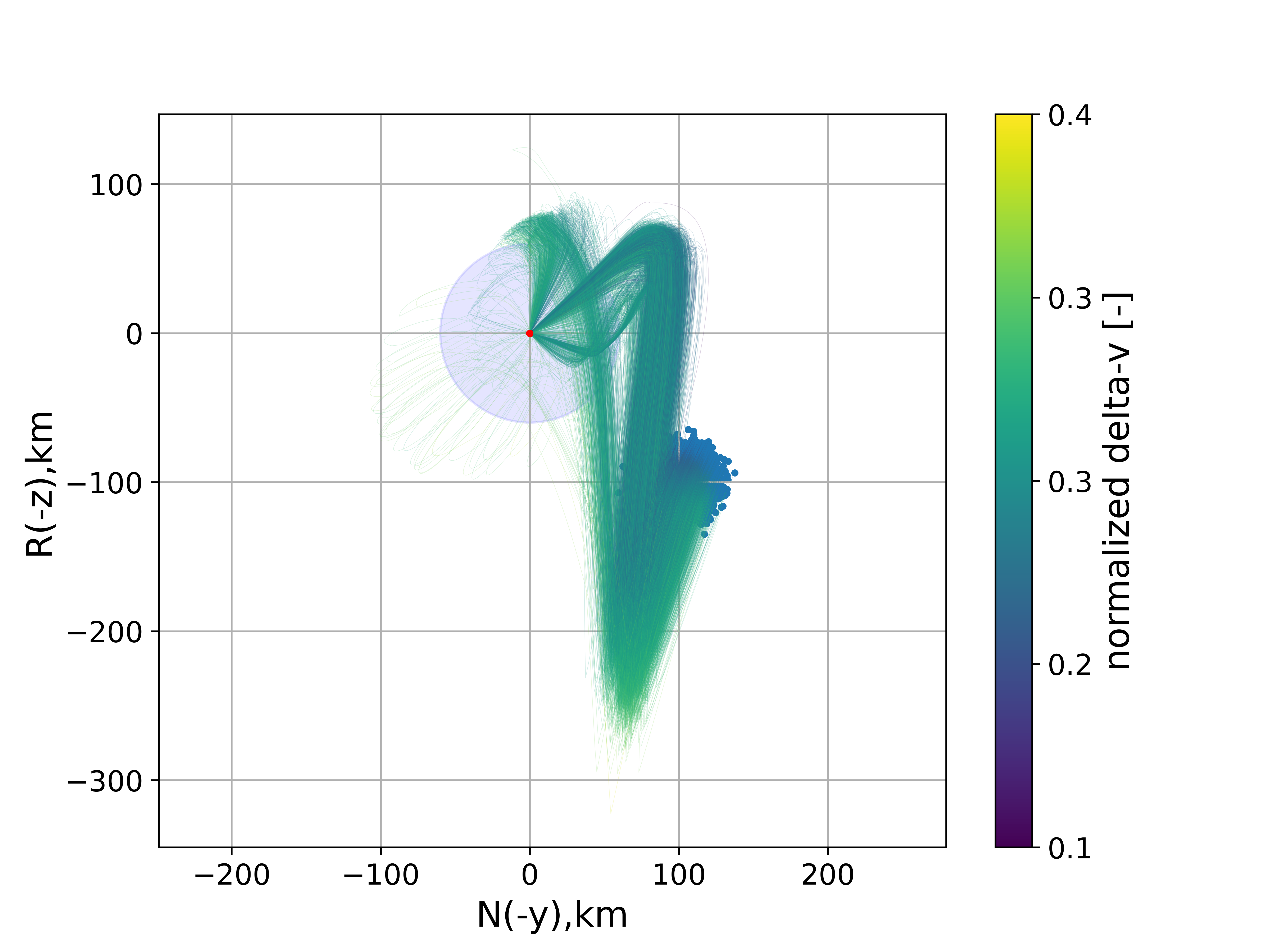}
        \caption{R-N Plane}
        \label{fig:sub2}
    \end{subfigure}
    
    \vspace{1em} 
    
    \begin{subfigure}[b]{0.48\textwidth}
        \centering
        \includegraphics[width=\linewidth]{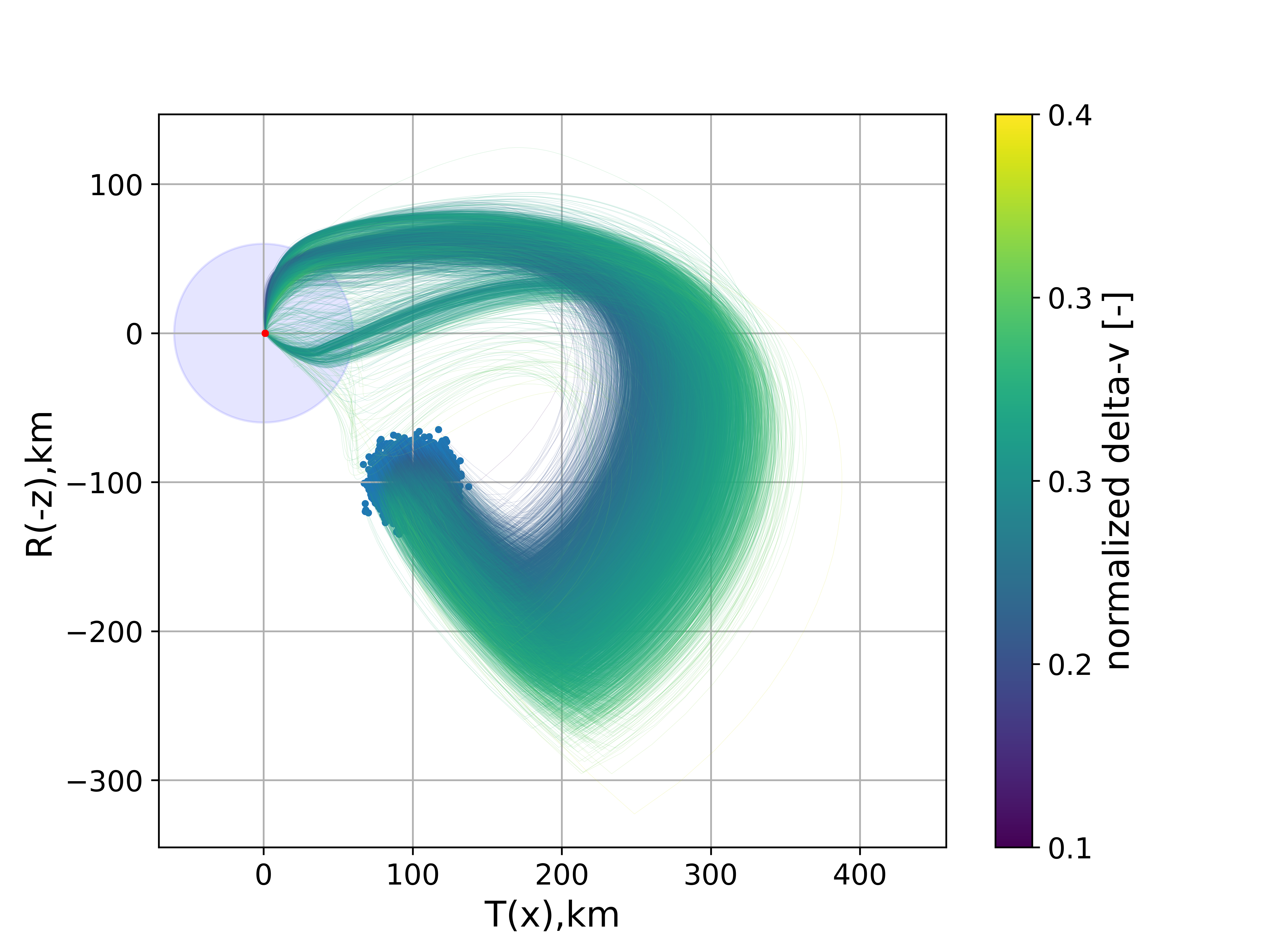}
        \caption{R-T Plane}
        \label{fig:sub3}
    \end{subfigure}
    \hfill
    \begin{subfigure}[b]{0.48\textwidth}
        \centering
        \includegraphics[width=\linewidth]{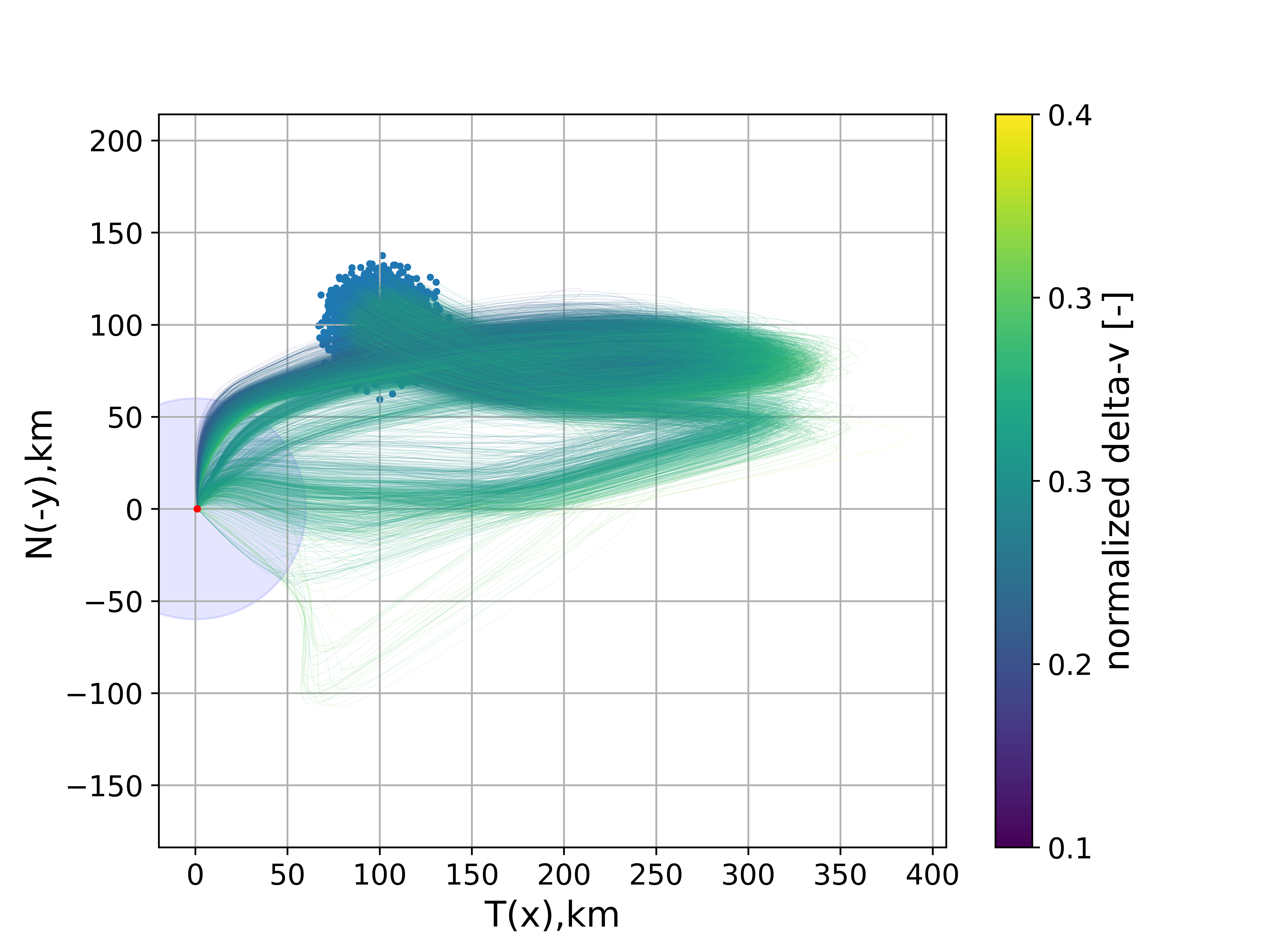}
        \caption{T-N Plane}
        \label{fig:sub4}
    \end{subfigure}
    
    \caption{Sampling of 10,000 trajectories out of the full dataset of 100,000.}
    \label{fig:dataset}
\end{figure}

\begin{figure}[htb!]
    \centering
    \includegraphics[width=0.7\linewidth]{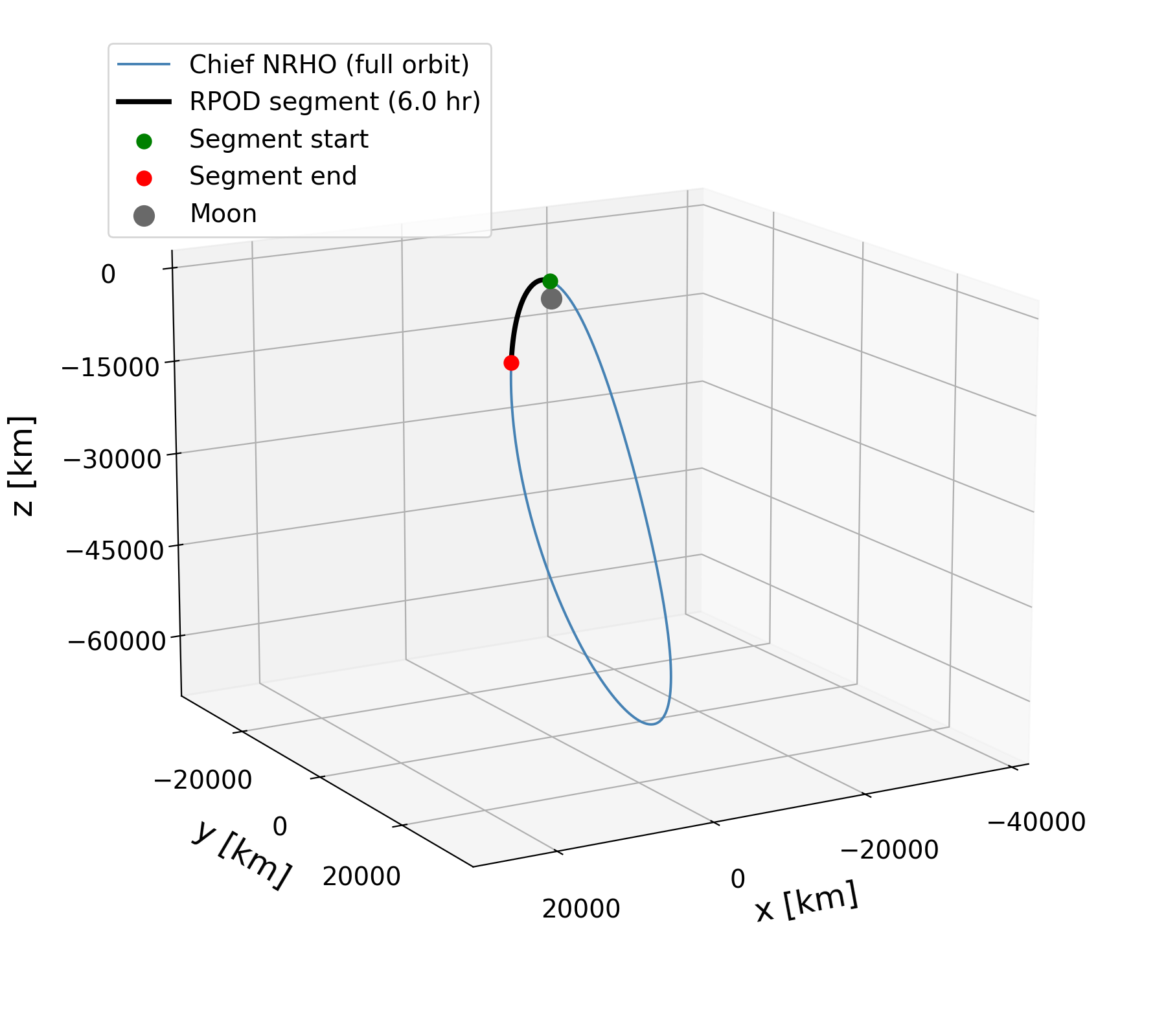}
    \caption{Chief NRHO with the segment used in this rendezvous scenario highlighted. Displayed in the Moon-centered synodic frame.}
    \label{fig:chieforbit}
\end{figure}

The first step in the evaluation process is to generate a dataset of solved SCP trajectories that encompasses all the desired variations of the cislunar rendezvous scenario, and then to partition the full dataset into training and testing sets.
The parameters of the dataset are given in Table \ref{tab:dataset}.
\begin{table}[htb!]
    \centering
    \begin{tabular}{c|c}
        Parameter &  Value\\
        \hline
        \hline
        $N_d$ [-] & 100,000 \\
        train-test split [\%] & 90 / 10\\
        \hline
        $\bm{r}$ [km] & [ 151.101,  -51.965,  2777.578] \\
        $\dot{\bm{r}}$ [km/s] & [-0.00123, -1.8499, -0.01776] \\
        $x_i$ [km] & $\mathcal{N}(100, 10)$ \\
        $y_i$ [km] & $\mathcal{N}(-100, 10)$ \\
        $z_i$ [km] & $\mathcal{N}(100, 10)$ \\
        $\dot{x}_i$ [m/s] & $\mathcal{N}(0, 0.05)$ \\
        $\dot{y}_i$ [m/s] & $\mathcal{N}(0, 0.05)$ \\
        $\dot{z}_i$ [m/s] & $\mathcal{N}(0,0.05)$ \\
        $\bm{x}_f$ [km; km/s] & [1, 0, 0, 0, 0, 0] \\
        \hline
    \end{tabular}
    \caption{Dataset distributions and presets for the defined cislunar rendezvous scenario. Absolute states ($\bm{r}$, $\dot{\bm{r}}$) are given in the Moon-centered synodic frame, and relative states ($\bm{x}_i = [x_i, y_i, z_i, \dot{x}_i, \dot{y}_i, \dot{z}_i]$, $\bm{x}_f$) are given in the chief-centered LVLH frame.}
    \label{tab:dataset}
\end{table}
Although the initial conditions of each sample in the dataset vary broadly, the constraints of the OCP that is solved are fixed.
The parameters of the OCP defined by Eq. \eqref{eq: OCP_nonlinear} are given by Table \ref{tab:ocp}. 
\begin{table}[htb!]
    \centering
    \begin{tabular}{c|c}
        Parameter & Value \\
        \hline
        \hline
        $r_{ae}$ [km] & 60 \\
        $r_{koz}$ [km] & 0.5 \\
        $N$ [-]  & 100 \\
        $n_{gate}$ [-]  & 80 \\
        $n_{safe}$ [-]  & 20 \\
        rendezvous horizon [hrs]  & 6 \\
        safety horizon [hrs] & 24  \\ 
        \hline
    \end{tabular}
    \caption{Preset OCP parameters for the defined cislunar rendezvous scenario.}
    \label{tab:ocp}
\end{table}
The rendezvous horizon is discretized into $N$ timesteps where impulsive control maneuvers can be applied, while the safety horizon is discretized into $n_{safe}$ timesteps.
The simulated dataset of trajectories is depicted in Figure \ref{fig:dataset}.
A 9:2 South NRHO is used as the reference orbit for the chief spacecraft. The mission segment along the chief orbit is selected intentionally to create a challenging rendezvous scenario. Figure \ref{fig:chieforbit} shows the chief orbit and the 6-hour segment where the rendezvous is performed.
The orbital radius with respect to the Moon is 2782 km at the start of the control window, and 18,160 km at the end of the control window, indicating major variation in dynamical regimes.
Optimization of the rendezvous trajectory and control inputs requires solving both the fast-moving dynamics at perilune and those further away from the Moon where the third-body perturbations from Earth start to have a more dominant effect.
This dataset is then used to train the ART-TWIN model so that it can generate warm-start trajectories autoregressively.

\subsection{ART-TWIN vs. ART-0 Evaluation}
From the test dataset, a subdataset of 100 trajectories was randomly selected to give quantifiable evidence of the advantage ART-TWIN provides over ART-0.

\begin{figure}[htb!]
    \centering
    \includegraphics[width=1\linewidth]{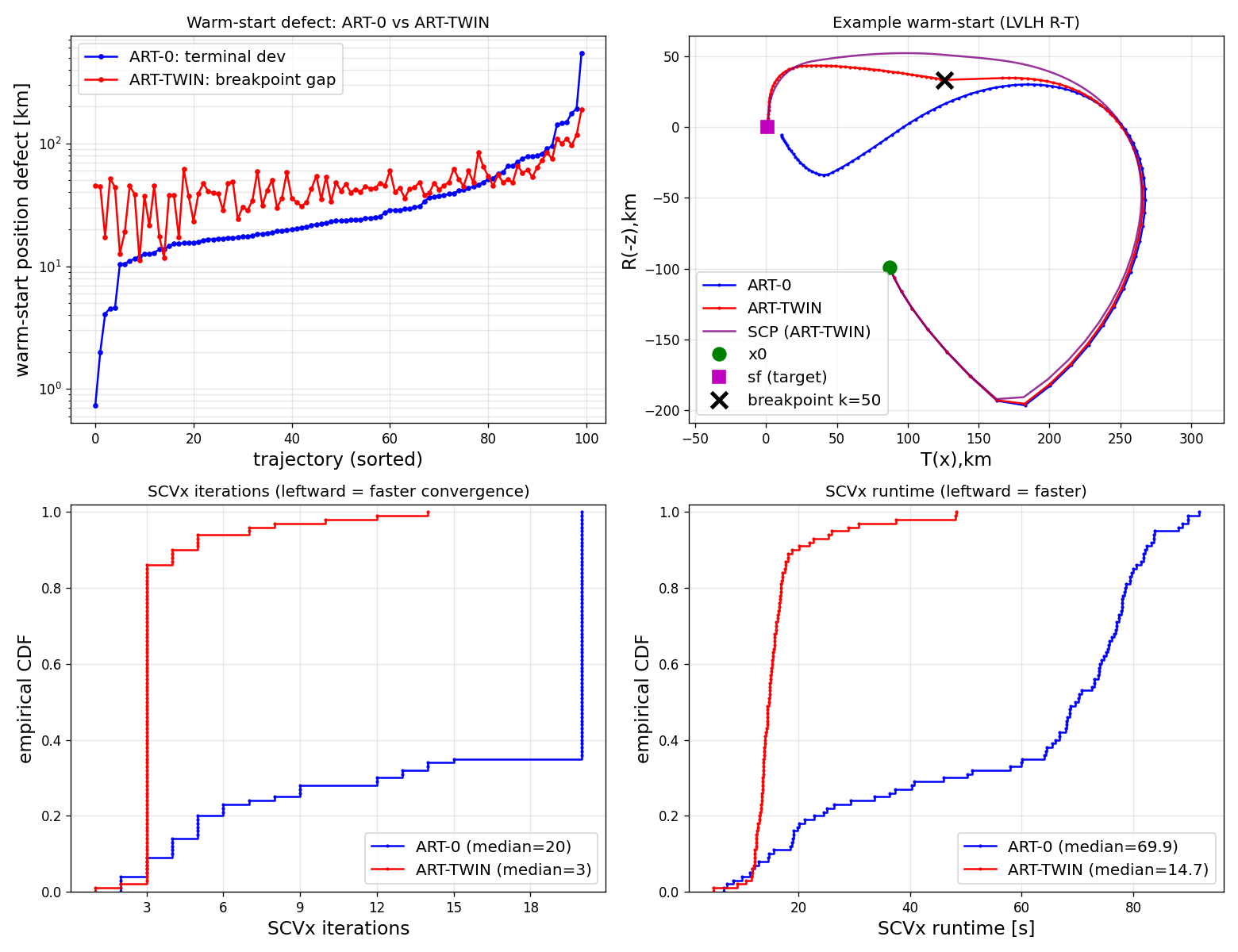}
    \caption{ART-TWIN vs. ART-0 comparison. Top left: Terminal state error of ART-0 warm-start compared to breakpoint error of ART-TWIN across the subdataset. Top right: Example of ART-0 and ART-TWIN warm-start trajectories, and fully solved SCP trajectory. Bottom left: cumulative distribution of SCP iterations conditioned on ART-TWIN and ART-0 warm-starts. Bottom right: cumulative distribution of SCP runtime conditioned on ART-TWIN and ART-0 warm-starts.}
    \label{fig:art_version_comp}
\end{figure}

Figure \ref{fig:art_version_comp} summarizes the results of this comparative analysis. The top right subplot shows the difference in open-loop trajectories generated by the two versions of ART for one selected case.
While there is a sharp jump at the midpoint to bridge the two arcs generated by ART-TWIN, the ART-0 trajectory has an excessive level of terminal state error.
It is important to make the distinction that the defect at the midpoint is only present in the ART-TWIN warm-start, and not in the final converged SCP trajectory.
As the same subfigure shows, the SCP trajectory warm-started by ART-TWIN is smooth and continuous, indicating that any breakpoint error was abated through SCP iterations.
This underscores the rationale behind ART-TWIN: Breakpoint state error is easier to resolve and correct than terminal state error.
To further investigate this claim, the top left subplot shows that for the majority of data, the warm-start defect is greater in ART-TWIN than ART-0.
Yet because that defect is at the midpoint rather than the terminal state, ART-TWIN remains more robust at providing an accurate warm start with reliable convergence than ART-0. 
Additionally, the defect in ART-TWIN is more predictable, generally staying on the order of 10 km, whereas the terminal state error in ART-0 varies widely across the subdataset, from as low as 1 km to over 100 km.
An example of the significant terminal state error exhibited by ART-0 can be clearly seen in the trajectory plot in the top right of Figure \ref{fig:art_version_comp}.
This reliability directly correlates to the efficiency performance metrics shown in the bottom row of subplots, where the full advantage of ART-TWIN is demonstrated.
About 90\% of the samples converge to an optimal solution via SCP in exactly 3 iterations when using ART-TWIN. However, when using ART-0, over 60\% of samples reach the maximum allowed number of iterations (20), meaning that they effectively fail to converge.
This underscores the motivation for enhancing ART via the addition of forward-backward shooting capabilities.
%



\subsection{ART-TWIN Warm-Starting Performance}

\begin{figure}[htb!]
    \centering
    \includegraphics[width=1\textwidth]{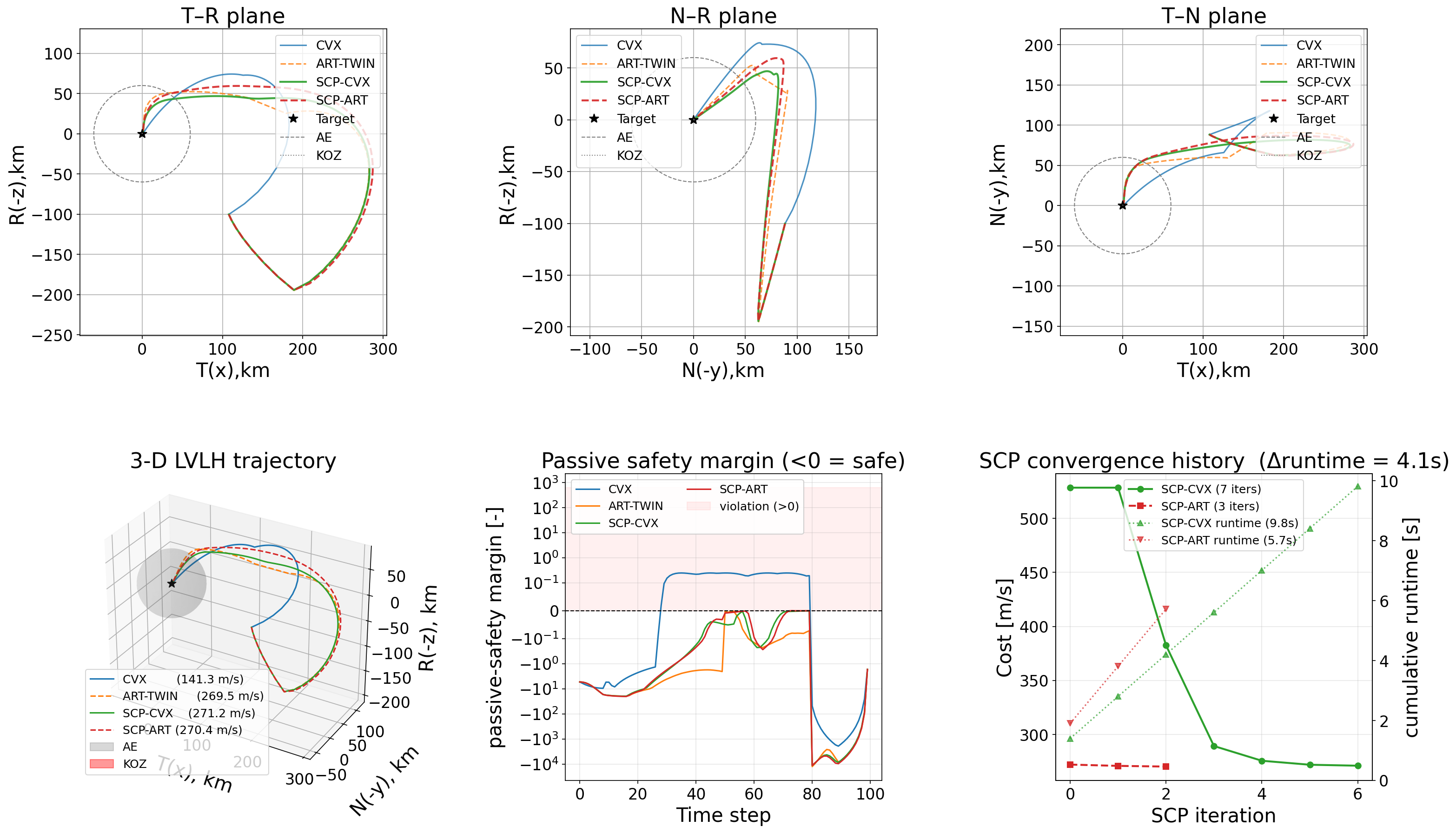}
    \caption{Comparison of CVX, ART, SCP-CVX, and SCP-ART solutions for one test sample. Top row: 2D trajectory plots in the LVLH frame. Bottom left: 3D trajectory plots in the LVLH frame. Bottom middle: Passive safety violations vs time. Bottom right: SCP iterations and runtime conditioned on the choice of warm start (CVX or ART-TWIN).}
    \label{fig:art_traj}
\end{figure}

ART-TWIN is then compared to the benchmark of a standard convex relaxation warm-start, obtained by solving Eq. \eqref{eq: OCP_linear}, to evaluate its performance.
As a remark on notation, the convex relaxation warm-start will be referred to going forward as ``CVX" and its associated SCP solution as ``SCP-CVX'', while the SCP solution warm-started by ART-TWIN will be referred to as ``SCP-ART". 
Figure \ref{fig:art_traj} shows a representative case from the cislunar rendezvous scenario test dataset.
In this case, the benefits of ART-TWIN as a warm-start to SCP are evident: The trajectory provided by ART-TWIN as an initial guess achieves fewer and less severe passive safety constraint violations, and the runtime of the ART-SCP solution is about 40\% faster than the CVX-SCP solution.
The midpoint discontinuity between the two arcs generated by ART-TWIN is evidently not a hindrance to rapid SCP convergence.
Qualitatively, it is clear why this is the case. 
Although the final SCP solutions are similar to each other, the initial guess generated by CVX takes a different path that violates the approach ellipsoid and keep-out-zone passive safety constraints, while the initial guess generated by ART-TWIN is already fairly close to the SCP solution. 
This also means that the delta-V estimate provided by ART-TWIN is much closer to the optimal solution, while the CVX estimate fails to provide a tight lower bound. 

\begin{figure}[htb!]
    \centering
    \includegraphics[width=1\linewidth]{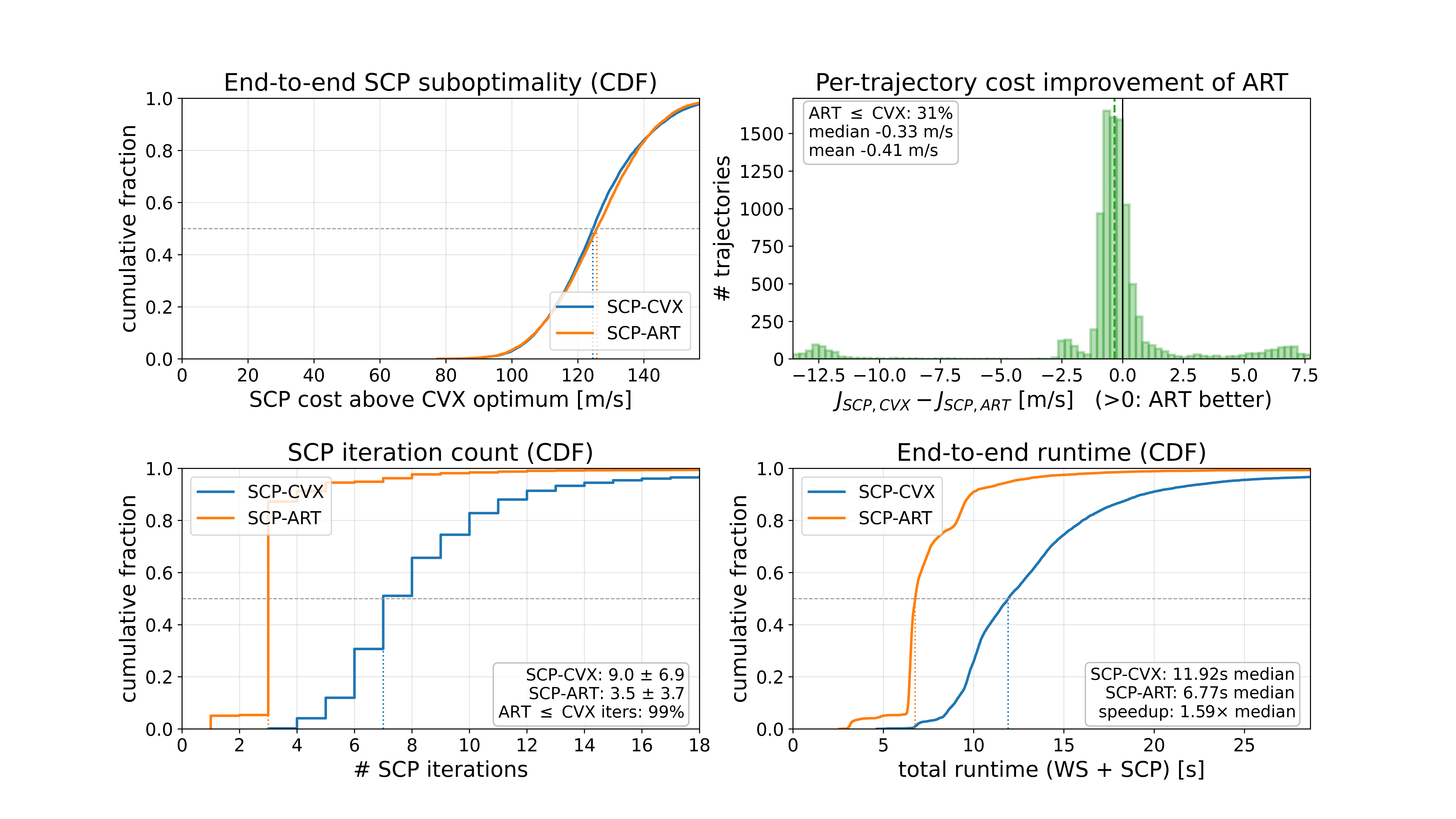}
    \caption{End-to-end ART-TWIN performance metrics, benchmarked against CVX, across the full test dataset.}
    \label{fig:art_metrics}
\end{figure}

An aggregated examination of the entire test dataset (10,000 held-out samples) indicates that this level of performance is a consistent trend. 
Figure \ref{fig:art_metrics} shows batched performance metrics and end-to-end results.
The greatest advantages can be seen in efficiency, as the subplots on the bottom row of Figure \ref{fig:art_metrics} illustrate that the ART-TWIN warm-start consistently outperforms the CVX warm-start in number of iterations and runtime, speeding up the SCP process by a median factor of 1.59x.
The top row displays cost and optimality metrics, which show small, but admittedly negligible, advantages for ART. 
The suboptimality gaps between the two methods are essentially the same, and the per-trajectory cost improvement from SCP-CVX to ART-SCP is only about 0.3-0.4 m/s on average, which is only a 1-2\% gain in total delta-V savings for the trajectories in this scenario.
Yet the benefits in efficiency provided by ART, specifically ART-TWIN, are clear. Note that in the efficiency diagnostics laid out in Figures \ref{fig:art_version_comp}, \ref{fig:art_traj}, and \ref{fig:art_metrics}, ART-TWIN actually outperforms ART-0 by more than it outperforms CVX.
In fact, an ART-0 warm-start would be incapable of matching the optimality and efficiency performance of CVX, due to its high likelihood of failed convergence.
Therefore, the architectural improvements that separate ART-TWIN from ART-0 are a necessary upgrade for ART to be considered a viable approach for cislunar trajectory optimization.

\section{Conclusion}

This work presents novel improvements to the Autonomous Rendezvous Transformer (ART), a transformer-boosted method of warm-starting sequential convex programming (SCP) to solve non-convex optimal control problems (OCPs). 
A forward-backward dual-arc generation enhancement is made to the original ART (ART-0) training and inference pipeline, enabling ART to handle OCPs that feature chaotic dynamical systems.
The enhanced version, ART-TWIN (Two-Way INference), is then successfully applied to cislunar relative motion, marking the first time ART has been extended to spacecraft rendezvous scenarios beyond Earth orbit.
ART-TWIN equals or outperforms both a convex relaxation and ART-0 as a warm start in every key efficiency and optimality metric, demonstrating a clear advantage.
Notably too, this performance is achieved in a rendezvous scenario that begins at perilune, a notoriously difficult region for controllability due to the rapidly varying dynamics present.
Furthermore, the promise ART-TWIN has shown in the cislunar dynamical regime demonstrates the value gain in using machine learning to augment conventional optimal control algorithms, and opens the door towards applying it to other chaotic and nonlinear systems, even those outside the field of spacecraft trajectory optimization \cite{takubo_transformer-based_2026}.

Further work is required to make ART universally deployable across the entire cislunar domain.
It will need to be trained and tested on a more diverse and expansive set of rendezvous scenarios and chief orbits\cite{celestini_generalizable_2025}.
Also, this work does not incorporate navigation uncertainties or take into account perturbations such as higher-resolution lunar gravity modeling, solar third-body gravitational effects, or solar radiation pressure.
Chance-constrained and model predictive control (MPC) approaches could be useful for addressing these respective extensions\cite{celestini_transformer-based_2024,takubo_towards_2025}.
Leveraging the natural manifold dynamics of the Circular Restricted Three-Body Problem (CR3BP) to bound quasi-periodic relative motion, thus making passive safety constraints easier to satisfy and guarantee, could also be a useful next step, depending on the desired use case \cite{takubo_safe_2026,elliott2022describing}.

Ultimately, the contributions of this work advance the onboard capabilities of autonomous distributed space systems (DSS), which will be vital for rendezvous, proximity operations, and docking (RPOD), as well as space domain awareness (SDA), in cislunar space. 

\section{Acknowledgments}

This work is supported in part by the National Science Foundation Graduate Research Fellowship Program under Grant No. DGE-2146755. This article solely reflects the opinions and conclusions of its authors and not any of its sponsors.
The authors thank Tommaso Guffanti, Daniele Gammelli, and other CAESAR researchers for their prior contributions in developing the original ART codebase.
AI tools, specifically Claude Sonnet 4.6 and Opus 4.8, were used to help debug and prototype the Python source code for this project. No AI tools were involved in the actual writing of this manuscript.

\section{Appendix}

\subsection{Transformer Model Hyperparameters}

The Autonomous Rendezvous Transformer (ART) presented in this
work is implemented in PyTorch \cite{paszke2017automatic} and builds on Hugging Face’s transformers library \cite{wolf-etal-2020-transformers}.
Table \ref{tab:art_params} presents the hyperparameter settings used in this work.
\begin{table}[ht]
\centering
\caption{ART hyperparameters.}
\begingroup
\begin{tabular}{l|l}
     Hyperparameter & Value \\
    \hline
    \hline
     Number of layers & 6\\
     Number of attention heads & 6 \\
    Embedding dimension & 384 \\
     Batch size& 4 \\
    Context length $K$ & 100 \\
    Nonlinearity & ReLU\\
    Dropout & 0.1\\
    Learning rate & 3e-5\\
    Grad norm clip & 1.0 \\
    Learning rate decay & None \\
    Gradient accumulation iters & 8 \\
    \hline
    \end{tabular}
    \label{tab:art_params}
    \endgroup
\end{table}
\subsection{SCP Hyperparameters}
Table \ref{tab:scp_params} presents the SCP hyperparameters used in this paper.
The nomenclature corresponds to the original SCvx* algorithm\cite{Oguri_2023}.
\begin{table}[ht]
\centering
\caption{SCvx* hyperparameters.}
\begingroup
\begin{tabular}{l|l}
     Hyperparameter & Value \\
    \hline
    \hline
    $\epsilon_{opt}$  & 1e-3\\
    $\epsilon_{feas}$ & 1e-3 \\
    $\rho_0$          & 0.0 \\
    $\rho_1$          & 0.25 \\
    $\rho_2$          & 0.7 \\
    $\alpha_1$        & 2 \\
    $\alpha_2$        & 2 \\
    $\beta$           & 1.5 \\
    $\gamma$          & 0.9 \\
    $r_\text{min}$    & 1e-6 \\
    $r_\text{max}$    & 10 \\
    $r^{(0)}$         & 0.5 \\
    $w^{(0)}$         & 10 \\
    max. iters.       & 20 \\
    \hline
    \end{tabular}
    \label{tab:scp_params}
    \endgroup
\end{table}

\appendix

\bibliographystyle{AAS_publication}   
\bibliography{references}   

@article{franzini_relative_2019,
	title = {Relative {Motion} {Dynamics} in the {Restricted} {Three}-{Body} {Problem}},
	volume = {56},
	issn = {0022-4650, 1533-6794},
	url = {https://arc.aiaa.org/doi/10.2514/1.A34390},
	doi = {10.2514/1.A34390},
	language = {en},
	number = {5},
	urldate = {2024-09-13},
	journal = {Journal of Spacecraft and Rockets},
	author = {Franzini, Giovanni and Innocenti, Mario},
	month = sep,
	year = {2019},
	pages = {1322--1337}
}

@article{vasile2003optimizing,
  title={Optimizing Low-Thrust and Gravity Assist Maneuvers to Design Interplanetary Trajectories},
  author={Vasile, Massimiliano and Bernelli-Zazzera, Franco},
  journal={The Journal of the Astronautical Sciences},
  volume={51},
  pages={13--35},
  year={2003},
  publisher={Springer}
}

@inproceedings{takubo2023optimization,
  title={Optimization Of Earth-Moon Low-Thrust-Enhanced Low-Energy Transfer},
  author={Takubo, Yuji and Shimane, Yuri and Ho, Koki},
  year = {2023},
  booktitle={2023 AAS/AIAA Astrodynamics Specialist Conference}
}

@article{sidhoum2024indirect,
  title={Indirect Forward--Backward Shooting for Low-Thrust Trajectory Optimization in Complex Dynamics},
  author={Sidhoum, Yanis and Oguri, Kenshiro},
  journal={Journal of Guidance, Control, and Dynamics},
  pages={1--9},
  year={2024},
  publisher={American Institute of Aeronautics and Astronautics}
}

@article{pierson1994three,
  title={Three-Stage Approach to Optimal Low-Thrust Earth-Moon Trajectories},
  author={Pierson, Bion L and Kluever, Craig A},
  journal={Journal of Guidance, Control, and Dynamics},
  volume={17},
  number={6},
  pages={1275--1282},
  year={1994}
}

@article{zhang2023stochastic,
  title={Stochastic Trajectory Optimization for 6-DOF Spacecraft Autonomous Rendezvous and Docking With Nonlinear Chance Constraints},
  author={Zhang, Yanquan and Cheng, Min and Nan, Bin and Li, Shunli},
  journal={Acta Astronautica},
  volume={208},
  pages={62--73},
  year={2023},
  publisher={Elsevier}
}

@article{guffanti2023passively,
  title={Passively Safe and Robust Multi-Agent Optimal Control With Application to Distributed Space Systems},
  author={Guffanti, Tommaso and D’Amico, Simone},
  journal={Journal of Guidance, Control, and Dynamics},
  volume={46},
  number={8},
  pages={1448--1469},
  year={2023},
  publisher={American Institute of Aeronautics and Astronautics}
}

@article{yamanaka_new_2002,
	title = {New {State} {Transition} {Matrix} for {Relative} {Motion} on an {Arbitrary} {Elliptical} {Orbit}},
	volume = {25},
	issn = {0731-5090, 1533-3884},
	url = {https://arc.aiaa.org/doi/10.2514/2.4875},
	doi = {10.2514/2.4875},
	language = {en},
	number = {1},
	urldate = {2024-09-15},
	journal = {Journal of Guidance, Control, and Dynamics},
	author = {Yamanaka, Koji and Ankersen, Finn},
	month = jan,
	year = {2002},
	pages = {60--66}
}

@article{clohessy_terminal_1960,
	title = {Terminal {Guidance} {System} for {Satellite} {Rendezvous}},
	volume = {27},
	issn = {1936-9999},
	url = {https://arc.aiaa.org/doi/10.2514/8.8704},
	doi = {10.2514/8.8704},
	language = {en},
	number = {9},
	urldate = {2024-09-16},
	journal = {Journal of the Aerospace Sciences},
	author = {Clohessy, W. H. and Wiltshire, R. S.},
	month = sep,
	year = {1960},
	pages = {653--658}
}

@inproceedings{banerjee_learning-based_2020,
	title = {Learning-Based {Warm}-{Starting} for {Fast} {Sequential} {Convex} {Programming} and {Trajectory} {Optimization}},
	url = {https://ieeexplore.ieee.org/document/9172293/?arnumber=9172293},
	doi = {10.1109/AERO47225.2020.9172293},
	urldate = {2024-09-17},
	booktitle = {2020 {IEEE} {Aerospace} {Conference}},
	author = {Banerjee, Somrita and Lew, Thomas and Bonalli, Riccardo and Alfaadhel, Abdulaziz and Alomar, Ibrahim Abdulaziz and Shageer, Hesham M and Pavone, Marco},
	month = mar,
	year = {2020},
	note = {ISSN: 1095-323X},
	pages = {1--8},
}

@inproceedings{malyuta2020fast,
  title={Fast Trajectory Optimization via Successive Convexification for Spacecraft Rendezvous With Integer Constraints},
  author={Malyuta, Danylo and Reynolds, Taylor and Szmuk, Michael and Acikmese, Behcet and Mesbahi, Mehran},
  booktitle={AIAA Scitech 2020 Forum},
  pages={0616},
  month = jan,
  year={2020}
}

@inproceedings{starek_real-time_2016,
	address = {Big Sky, MT, USA},
	title = {Real-Time, Propellant-Optimized Spacecraft Motion Planning Under {Clohessy}-{Wiltshire}-{Hill} dynamics},
	isbn = {978-1-4673-7676-1},
	url = {http://ieeexplore.ieee.org/document/7500704/},
	doi = {10.1109/AERO.2016.7500704},
	language = {en},
	urldate = {2024-09-17},
	booktitle = {2016 {IEEE} {Aerospace} {Conference}},
	publisher = {IEEE},
	author = {Starek, Joseph A. and Schmerling, Edward and Maher, Gabriel D. and Barbee, Brent W. and Pavone, Marco},
	month = mar,
	year = {2016},
	pages = {1--16},
}

@article{koenig2017new,
  title={New State Transition Matrices for Spacecraft Relative Motion in Perturbed Orbits},
  author={Koenig, Adam W and Guffanti, Tommaso and D’Amico, Simone},
  journal={Journal of Guidance, Control, and Dynamics},
  volume={40},
  number={7},
  pages={1749--1768},
  year={2017},
  publisher={American Institute of Aeronautics and Astronautics}
}

@article{elliott2022describing,
  title={Describing Relative Motion Near Periodic Orbits via Local Toroidal Coordinates},
  author={Elliott, Ian and Bosanac, Natasha},
  journal={Celestial Mechanics and Dynamical Astronomy},
  volume={134},
  number={2},
  pages={19},
  year={2022},
  publisher={Springer}
}

@inproceedings{khoury_relative_2022,
	title = {Relative {Guidance}, {Navigation}, and {Control}, in {Multibody} {Gravitational} {Regimes}},
    booktitle = {2022 AAS/AIAA Astrodynamics Specialist Conference},
	language = {en},
	author = {Khoury, Fouad and Lippe, Corinne},
	year = {2022},
    month = aug,
}

@inproceedings{khoury_relative_2024,
	address = {Broomfield, CO},
	title = {Relative-{Navigation}-{Informed} {Formation} {Configuration} and {Station}-{Keeping} {Strategies} in the {Earth}-{Moon}-{Sun} {System}},
	language = {en},
	booktitle = {2024 {AAS}/{AIAA} {Astrodynamics} {Specialist} {Conference}},
	author = {Khoury, Fouad and Power, Rolfe and Ghosh, Pradipto and Romero, Juan Ojeda},
	year = {2024},
}

@article{vela2026application,
  title={Application of Fundamental Modal Solutions to Relative Dynamics in the Cislunar Environment},
  author={Vela, Claudio and Opromolla, Roberto and Fasano, Giancarmine and Schaub, Hanspeter},
  journal={Journal of Guidance, Control, and Dynamics},
  volume={49},
  number={2},
  pages={344--358},
  year={2026},
  publisher={American Institute of Aeronautics and Astronautics}
}

@article{breger2008safe,
  title={Safe Trajectories for Autonomous Rendezvous Of Spacecraft},
  author={Breger, Louis and How, Jonathan P},
  journal={Journal of Guidance, Control, and Dynamics},
  volume={31},
  number={5},
  pages={1478--1489},
  year={2008}, 
  doi={10.2514/1.29590}
}

@article{howell1984three,
  title={Three-Dimensional Periodic Halo Orbits},
  author={Howell, Kathleen Connor},
  journal={Celestial mechanics},
  volume={32},
  pages={53},
  year={1984}
}

@mastersthesis{zimovan2017characteristics,
  title={Characteristics and Design Strategies for Near Rectilinear Halo Orbits Within the Earth-Moon System},
  author={Zimovan, Emily M},
  year={2017},
  school={Purdue University}
}

@article{takubo_safe_2026,
author = {Takubo, Yuji and Manuel, Walter and Foss, Ethan and D’Amico, Simone},
title = {Safe and Optimal N-Spacecraft Swarm Reconfiguration in Non-Keplerian Cislunar Orbits},
journal = {Journal of Guidance, Control, and Dynamics},
volume = {49},
number = {2},
pages = {305-324},
year = {2026},
doi = {10.2514/1.G009290},
URL = {https://doi.org/10.2514/1.G009290},
eprint = {https://doi.org/10.2514/1.G009290}
}

@inproceedings{foss_long-duration_2025,
    author = {Foss, Ethan and Takubo, Yuji and D'Amico, Simone},
    title = {Long-{Duration} {Station}-{Keeping} {Strategy} for {Cislunar} {Spacecraft} {Formations}},
    booktitle = {2025 AAS/AIAA Astrodynamics Specialist Conference},
    address = {Boston, MA},
    month = aug,
    year = {2025}
}

@inproceedings{hunter_optimal_2025,
	title = {Optimal {Impulsive} {Control} of {Cislunar} {Relative} {Motion} using {Reachable} {Set} {Theory}},
	author = {Hunter, Matthew and Manuel, Walter J. and D'Amico, Simone},
    booktitle = {2025 AAS/AIAA Astrodynamics Specialist Conference},
    address = {Boston, MA},
    month = aug,
    year = {2025}
}

@inproceedings{paszke2017automatic,
  title={Automatic Differentiation in PyTorch},
  author={Paszke, Adam and Gross, Sam and Chintala, Soumith and Chanan, Gregory and Yang, Edward and DeVito, Zachary and Lin, Zeming and Desmaison, Alban and Antiga, Luca and Lerer, Adam},
  booktitle={NIPS-W},
  year={2017}
}

@inproceedings{wolf-etal-2020-transformers,
    title = "Transformers: State-Of-The-Art Natural Language Processing",
    author = "Thomas Wolf and Lysandre Debut and Victor Sanh and Julien Chaumond and Clement Delangue and Anthony Moi and Pierric Cistac and Tim Rault and Rémi Louf and Morgan Funtowicz and Joe Davison and Sam Shleifer and Patrick von Platen and Clara Ma and Yacine Jernite and Julien Plu and Canwen Xu and Teven Le Scao and Sylvain Gugger and Mariama Drame and Quentin Lhoest and Alexander M. Rush",
    booktitle = "Proceedings Of the 2020 Conference On Empirical Methods in Natural Language Processing: System Demonstrations",
    month = oct,
    year = "2020",
    address = "Online",
    publisher = "Association for Computational Linguistics",
    url = "https://aclanthology.org/2020.emnlp-demos.6/",
    pages = "38--45"
}

@inproceedings{Oguri_2023,
   title={Successive Convexification With Feasibility Guarantee via Augmented Lagrangian for Non-Convex Optimal Control Problems},
   url={http://dx.doi.org/10.1109/CDC49753.2023.10383462},
   DOI={10.1109/cdc49753.2023.10383462},
   booktitle={2023 62nd IEEE Conference On Decision and Control (CDC)},
   publisher={IEEE},
   author={Oguri, Kenshiro},
   year={2023},
   month=Dec, pages={3296–3302} }

@misc{takubo_transformer-based_2026,
    title = {Transformer-{Based} {Warm}-{Starting} for {Feasible} and {Optimal} {Terminal} {Approach} to {Tumbling} {Objects} with {Space} {Manipulators}},
    url = {http://arxiv.org/abs/2606.17317},
    doi = {10.48550/arXiv.2606.17317},
    urldate = {2026-06-22},
    publisher = {arXiv},
    author = {Takubo, Yuji and Adang, Maximilian and Schwager, Mac and D'Amico, Simone},
    month = jun,
    year = {2026},
    note = {arXiv:2606.17317 [cs.RO]},
}

@inproceedings{takubo_agile_2026,
    title = {Agile {Tradespace} {Exploration} for {Space} {Rendezvous} {Mission} {Design} via {Transformers}},
    issn = {2996-2358},
    url = {https://ieeexplore.ieee.org/document/11519905/},
    doi = {10.1109/AERO66936.2026.11519905},
    urldate = {2026-07-03},
    booktitle = {2026 {IEEE} {Aerospace} {Conference}},
    author = {Takubo, Yuji and Gammelli, Daniele and Pavone, Marco and D'Amico, Simone},
    month = mar,
    year = {2026},
    note = {ISSN: 2996-2358},
    pages = {1--14},
}

@inproceedings{takubo_towards_2025,
    title = {Towards {Robust} {Spacecraft} {Trajectory} {Optimization} via {Transformers}},
    issn = {2996-2358},
    url = {https://ieeexplore.ieee.org/document/11068406/},
    doi = {10.1109/AERO63441.2025.11068406},
    urldate = {2026-07-03},
    booktitle = {2025 {IEEE} {Aerospace} {Conference}},
    author = {Takubo, Yuji and Guffanti, Tommaso and Gammelli, Daniele and Pavone, Marco and D'Amico, Simone},
    month = mar,
    year = {2025},
    note = {ISSN: 2996-2358},
    pages = {1--13},
}

@article{celestini_transformer-based_2024,
    title = {Transformer-{Based} {Model} {Predictive} {Control}: {Trajectory} {Optimization} via {Sequence} {Modeling}},
    volume = {9},
    issn = {2377-3766},
    shorttitle = {Transformer-{Based} {Model} {Predictive} {Control}},
    url = {https://ieeexplore.ieee.org/document/10685132/},
    doi = {10.1109/LRA.2024.3466069},
    number = {11},
    urldate = {2026-07-03},
    journal = {IEEE Robotics and Automation Letters},
    author = {Celestini, Davide and Gammelli, Daniele and Guffanti, Tommaso and D'Amico, Simone and Capello, Elisa and Pavone, Marco},
    month = nov,
    year = {2024},
    pages = {9820--9827},
}

@inproceedings{celestini_generalizable_2025,
    title = {Generalizable {Spacecraft} {Trajectory} {Generation} via {Multimodal} {Learning} with {Transformers}},
    issn = {2378-5861},
    url = {https://ieeexplore.ieee.org/document/11108053/},
    doi = {10.23919/ACC63710.2025.11108053},
    urldate = {2026-07-03},
    booktitle = {2025 {American} {Control} {Conference} ({ACC})},
    author = {Celestini, Davide and Afsharrad, Amirhossein and Gammelli, Daniele and Guffanti, Tommaso and Zardini, Gioele and Lall, Sanjay and Capello, Elisa and D’Amico, Simone and Pavone, Marco},
    month = jul,
    year = {2025},
    note = {ISSN: 2378-5861},
    pages = {3558--3565},
}

@article{briden_constraint-informed_2025,
    title = {Constraint-{Informed} {Learning} for {Warm}-{Starting} {Trajectory} {Optimization}},
    volume = {48},
    issn = {0731-5090, 1533-3884},
    url = {https://arc.aiaa.org/doi/10.2514/1.G008791},
    doi = {10.2514/1.G008791},
    language = {en},
    number = {10},
    urldate = {2026-07-03},
    journal = {Journal of Guidance, Control, and Dynamics},
    author = {Briden, Julia and Choi, Changrak and Yun, Kyongsik and Linares, Richard and Cauligi, Abhishek},
    month = oct,
    year = {2025},
    pages = {2272--2287},
}

@inproceedings{guffanti_transformers_2024,
    title = {Transformers for {Trajectory} {Optimization} with {Application} to {Spacecraft} {Rendezvous}},
    issn = {1095-323X},
    url = {https://ieeexplore.ieee.org/document/10521334/},
    doi = {10.1109/AERO58975.2024.10521334},
    urldate = {2026-07-03},
    booktitle = {2024 {IEEE} {Aerospace} {Conference}},
    author = {Guffanti, Tommaso and Gammelli, Daniele and D’Amico, Simone and Pavone, Marco},
    month = mar,
    year = {2024},
    note = {ISSN: 1095-323X},
    pages = {1--13},
}

@inproceedings{vaswani_attention_2017,
    title = {Attention is {All} you {Need}},
    volume = {30},
    url = {https://proceedings.neurips.cc/paper/2017/hash/3f5ee243547dee91fbd053c1c4a845aa-Abstract.html},
    urldate = {2026-07-03},
    booktitle = {Advances in {Neural} {Information} {Processing} {Systems}},
    publisher = {Curran Associates, Inc.},
    author = {Vaswani, Ashish and Shazeer, Noam and Parmar, Niki and Uszkoreit, Jakob and Jones, Llion and Gomez, Aidan N and Kaiser, {\L} ukasz and Polosukhin, Illia},
    year = {2017},
}

@phdthesis{malyuta_thesis_2021,
    address = {United States -- Washington},
    type = {Ph.{D}. thesis},
    title = {Convex {Optimization} in a {Nonconvex} {World}: {Applications} for {Aerospace} {Systems}},
    copyright = {Database copyright ProQuest LLC; ProQuest does not claim copyright in the individual underlying works.},
    isbn = {979-8-5355-0153-6},
    shorttitle = {Convex {Optimization} in a {Nonconvex} {World}},
    url = {https://www.proquest.com/docview/2568283655/abstract/A241A5E24DE94870PQ/1},
    language = {English},
    urldate = {2026-07-06},
    school = {University of Washington},
    author = {Malyuta, Danylo},
    year = {2021},
}

@inproceedings{guzzetti_stationkeeping_2017,
    title = {{Stationkeeping} {Analysis} {for} {Spacecraft} {in} {Lunar} {Near} {Rectilinear} {Halo} {Orbits}},
    language = {en},
    booktitle = {2017 AAS/AIAA Space Flight Mechanics Meeting},
    month = jan,
    year = {2017},
    author = {Guzzetti, Davide and Zimovan, Emily M and Howell, Kathleen C and Davis, Diane C},
}

@misc{shimane_autonomous_2026,
    title = {Autonomous {Navigation} and {Station}-{Keeping} on {Near}-{Rectilinear} {Halo} {Orbits}},
    url = {http://arxiv.org/abs/2512.01182},
    doi = {10.48550/arXiv.2512.01182},
    urldate = {2026-07-07},
    publisher = {arXiv},
    author = {Shimane, Yuri and Berntorp, Karl and Cairano, Stefano Di and Weiss, Avishai},
    month = may,
    year = {2026},
    note = {arXiv:2512.01182 [eess.SY]},
}

@article{black_fragmentation_2025,
    title = {Fragmentation characterization in the circular restricted three body problem for cislunar space domain awareness},
    volume = {75},
    issn = {0273-1177},
    url = {https://www.sciencedirect.com/science/article/pii/S0273117724009189},
    doi = {10.1016/j.asr.2024.08.076},
    number = {1},
    urldate = {2026-07-07},
    journal = {Advances in Space Research},
    author = {Black, Arly and Frueh, Carolin},
    month = jan,
    year = {2025},
    pages = {1177--1204},
}

@misc{liyanaarachchi_6g_2024,
    title = {{6G} at $\frac{1}{6}g$: {The} {Future} of {Cislunar} {Communications}},
    shorttitle = {{6G} at $\frac{1}{6}g$},
    url = {http://arxiv.org/abs/2407.16672},
    doi = {10.48550/arXiv.2407.16672},
    urldate = {2026-07-07},
    publisher = {arXiv},
    author = {Liyanaarachchi, Sahan and Mitrolaris, Stavros and Mitra, Purbesh and Ulukus, Sennur},
    month = jul,
    year = {2024},
    note = {arXiv:2407.16672 [cs.IT]},
}

@inproceedings{sidhoum_low_2023,
author = {Sidhoum, Yanis and Oguri, Kenshiro},
year = {2023},
month = aug,
booktitle = {2023 AAS/AIAA Astrodynamics Specialist Conference},
title = {Low-thrust Trajectory Optimization for Enceladus Exploration using Indirect Forward-Backward Shooting}
}

\end{document}